\RequirePackage{fix-cm}
\documentclass[smallextended]{svjour3}       
\smartqed  
\usepackage{graphicx}

\RequirePackage[OT1]{fontenc} \RequirePackage[colorlinks]{hyperref}
\RequirePackage{hypernat}

\usepackage{latexsym,enumerate}
\usepackage{amsmath,amsopn,amstext,amscd,amsfonts,amssymb,
pdfsync,amsbsy}
\usepackage{mathrsfs}
\usepackage{fullpage}
\usepackage{graphicx}
\usepackage[usenames]{color}

\newtheorem{rem}{Remark}

\let\epsilon=\varepsilon
\let\eps=\epsilon
\let\phi=\varphi
\let\al=\alpha

\let\la=\lambda \let\La=\Lambda 
\let\tilde=\widetilde

\newcommand{\field}[1]{\mathbb{#1}}
\newcommand{\R}{\field{R}}

\newcommand{\N}{\field{N}}
\newcommand{\Z}{\field{Z}}

\newcommand{\EE}{\field{E}}

\renewcommand{\L}{\field{L}}
\newcommand{\bH}{{\field H}}
\newcommand{\bS}{{\field S}}

\newcommand{\F}{{\mathscr{F}}}

\newcommand{\A}{{\mathcal A}}
\newcommand{\B}{{\mathcal B}}
\newcommand{\V}{{\mathcal V}}

\newcommand{\M}{{\mathcal M}}
\newcommand{\D}{\mathcal{D}}

\newcommand{\cH}{\mathcal{H}}

\newcommand{\cN}{{\mathcal N}}
\newcommand{\cT}{\mathcal{T}}

\newcommand{\cP}{{\mathcal P}}

\newcommand{\beqn}{\begin{equation}}
\newcommand{\eeqn}{\end{equation}}
\newcommand\eref[1]{(\ref{#1})}

\newcommand{\foral}{\forall \;}

\def\E{{\field E}}

\def\P{{\field P}}

\def\cE{{\field E}}

\def\lam{\lambda}

\def\C{\mathcal C}

\def\supp{{\rm{supp}\, }}
\def\UCB{{\rm UCB}}

\newcommand{\mb}{\field{B}}          
\newcommand{\enb}{\bar{\eps}_n}  
\newcommand{\rn}{\sqrt{n}}
\newcommand{\leqa}{\lesssim}
\newcommand{\pli}{+\infty}
\mathchardef\given="626A

\newcommand{\cX}{{\mathcal X}}

\begin{document}

\title{Thomas Bayes' walk on manifolds
\thanks{I.Castillo's work is partly supported by ANR Grant `Banhdits' ANR-2010-BLAN-0113-03.}
}


\author{Isma\"el Castillo        \and
        G\'erard Kerkyacharian \and
        Dominique Picard 
}


\institute{Isma\"el Castillo \at
              CNRS-LPMA Universities Paris 6 and 7,
              \email{ismael.castillo@upmc.fr}           
           \and
             G\'erard Kerkyacharian and Dominique Picard \at
             Universit\'e Paris Diderot - Paris 7, LPMA
             \email{kerk@math.jussieu.fr, picard@math.jussieu.fr} 
}

\date{Received: date / Accepted: date}

\maketitle

\begin{abstract}
Convergence of the Bayes posterior measure is considered in canonical statistical settings where observations sit on a geometrical object such as a compact manifold, or more generally on a compact metric space verifying some conditions. 
A natural geometric prior based on randomly rescaled solutions of the heat equation is considered. Upper and lower bound posterior contraction rates are derived.
\keywords{Bayesian nonparametrics \and Gaussian process priors \and Heat kernel}
\subclass{MSC 62G05 \and MSC 62G20}
\end{abstract}

\baselineskip=18 pt

\section{Introduction}

Let $\M$ be a compact metric space, equipped with a Borel measure $\mu$ and the corresponding Borel-sigma field. Let $\L^p:=\L^p(\M,\mu)$, $p\ge 1$ denote the space of $p$-integrable real functions defined on $\M$ with respect to $\mu$.

In this paper we investigate rates of contraction of posterior distributions for nonparametric
 models on geometrical structures such as 
\begin{enumerate}
\item Gaussian white noise on a compact metric space $\M$, where, for $n\ge 1$,
one observes 
\begin{equation*}
dX^{(n)}(x) = f(x)dx + \frac{1}{\rn} dZ(x),\quad x\in \M, 
\end{equation*}
where $f$ is in $\L^2$ and $Z$ is a white noise on $\M$.
\item Fixed design regression
where one observes, for $n\ge 1$,
\begin{equation*}
Y_i = f(x_i) + \eps_i,\quad 1\le i\le n.
\end{equation*}
The design points $\{x_i\}$ are fixed on $\M$ and the variables $\{\eps_i\}$
are assumed to be independent standard normal. 
\item Density estimation on a manifold
where the observations are a sample 
\begin{equation*}
(X_i)_{1\le i\le n}\ \sim\ f,
\end{equation*}
$X_1,\ldots,X_n$ are independent identically distributed $\M$-valued random variables with positive density function $f$ on $\M$. 
\end{enumerate}


Although an impressive amount of work has been done using frequentist approachs to estimation on manifolds, see \cite{MardiaJupp} and the references therein, we focus in this paper on the  Bayes posterior measure.
Works devoted to deeply understanding the behaviour of Bayesian nonparametric methods have recently experienced a considerable development in particular after the seminal works of A. W. van der Vaart, H. van Zanten, S. Ghosal and J. K. Ghosh \cite{ggv00}, \cite{shenwasserman01}. 
Especially, the class of Gaussian processes forms an important family of nonparametric prior distributions, for which precise rates have been obtained  \cite{vvvz}, see also \cite{ic08} for lower bound counterparts. In \cite{vvvz09}, the authors obtained adaptive performance up to logarithmic terms by introducing a random rescaling of a very smooth Gaussian random field. These results have been obtained on $[0,1]^d, d\ge 1$. Our point in this paper 
is to develop a Bayesian procedure adapted to the geometrical structure of the data.
Among the examples covered by our results, we can cite directional data corresponding to the spherical case and more generally the case of data supported by a compact manifold. 
 
We follow the illuminating approach of \cite{vvvz} and \cite{vvvz09} and use a fixed prior distribution, constructed by rescaling a smooth Gaussian random field. In our more general  setting, we show how the rescaling is made possible by introducing a notion of time decoupled from the underlying space. Another important difference brought by the geometrical nature of the problem is the underlying Gaussian process, which now originates from an harmonic analysis of the data space $\M$, with the rescaling naturally acting on the frequency domain.

We suppose that $\M$ is equipped with a positive self-adjoint operator $L$ such that the associated semi-group $e^{-tL}$, $t>0$,  the {\it heat kernel}, allows a  smooth functional calculus,  which in turn allows the construction of the Gaussian random field. Our prior can then be interpreted as a randomly rescaled (random) solution of the heat equation. 

We also took inspiration on earlier work by  \cite{angers-kim}, where the authors consider a  symmetry-adaptive Bayesian estimator in a regression framework. Precise minimax rates in the  $\L^2$-norm over Sobolev spaces of functions on compact connected orientable manifolds without boundary are obtained in \cite{efromovich}. We also mention a recent development by \cite{bdunson10}, where Bayesian consistency properties are derived for priors based on mixture of kernels over a compact manifold.

Here is an outline of the paper. 
We first detail in Section \ref{geo} the properties assumed on the structure $\M$ and the associated heat kernel allowing our construction and  give examples.
We then construct the associated Gaussian prior in Section \ref{section_rkhs} and prove some approximation and concentration properties typically needed to obtain rates of contraction of posterior distributions in Sections \ref{sec-prior} and \ref{section_gc}.
We then prove in Section \ref{section_rates} upper-bound rates in the three statistical examples detailed above. Lower bounds are considered in Section \ref{sec-lb}. 
The Appendices in Sections \ref{appendix-besov}, \ref{appendix-entropy} and \ref{appendix-manifold} contain respectively the definition of Besov spaces, the proofs of entropy results and a property of measure of balls on compact Riemannian manifolds.

The notation $\leqa$ means less than or equal to up to some universal constant. For any sequences of reals $(a_n)_{n\ge 0}$ and $(b_n)_{n\ge 0}$, the notation $a_n\sim b_n$  means that the sequences verify $c\le \liminf_{n} (b_n/a_n) \le \limsup_n (b_n/a_n) \le d$ for some positive constants $c,d$, and $a_n\ll b_n$ stands for $\lim_n (b_n/a_n)=0$. For any reals $a,b$, we denote $\min(a,b)=a\wedge b$ and $\max(a,b)=a\vee b$.

\section{The geometrical framework}\label{geo}

The squared-exponential covariance kernel introduced in \cite{vvvz09}, which gives rise to a particular Reproducing Kernel Hilbert Space (RKHS) see \cite{rkhs} and Section \ref{section_rkhs} below, has in fact a natural extension to more general metric spaces.  

Suppose $\L^2 = \oplus_{k\ge 0} \cH_{k}$, where the  $\cH _{k}$ are supposed to be  finite-dimensional subspaces of $\L^2$ consisting of continuous functions on $\M$, and orthogonal in $\L^2$. Then, the projector $P_{k}$ on $\cH_k$ is actually a kernel operator $P_k(x,y):=\sum_{1\le i \le dim(\cH_k)} e_k^i(x) e_k^i(y)$, where $\{e_k^i\}$ is any orthonormal basis of $\cH_k$; so it is obviously a positive-definite kernel. Also, given $\phi:\mathbb{N}\to (0,+\infty)$ and under a uniform convergence assumption, $K_\phi(x,y)=\sum_{k\ge 0}\phi(k)P_k(x,y)$ is a positive definite kernel which is the covariance kernel of  a Gaussian process.

Here, we will focus on the case where the subspaces $\cH_k=:\cH_{\la_k}$ are the eigenspaces of a 
self-adjoint positive operator $L$ and $\phi(k)=e^{-\la_k t}$, $t>0$, so that 
$\sum_{k} e^{-\la_k t} P_k(x,y)$ is actually the kernel  of the associated semi-group 
$e^{-tL}$.

This construction, under the following appropriate conditions, yields a natural generalisation of the squared-exponential covariance kernel on the real line.

\subsection{Compact metric doubling space} \label{sec-metspace}
 The open 
 balls of radius $r$  centered in $x \in \M$ are denoted by $B(x,r)$ and to simplify the notation we put
 $\mu(B(x,r)) =: |B(x,r)|.$ The metric is denoted by $\rho.$ For simplicity,  we can  impose, in the abstract proofs that $\mu(\M)=1= diam(\M)$.
 But of course, this is not the case in practical situations: the metric and the measure on a Riemannian compact manifold is not 
 normalized (see below 2.5).
We assume that $\M$ has the so called doubling property: i.e. 
there exists a constant $0<D <\infty$  {such that:}
 \begin{equation}\label{VD}
  \hbox{ for all }\; x \in \M,\; 0<r,  \quad 0< |B(x,2r)| \leq 2^D |B(x,r)|
 \end{equation}
\begin{rem}
As a simple consequence of \eref{VD} we have :
\begin{equation}
\hbox{ for all }\; x,y \in \M, \hbox{ for all }\; 0<r \leq R, \quad |B(x,R)| \leq ( \frac{2R}{r})^D (1+ \frac{\rho(x,y)}R)^D |B(y,r)|
\end{equation}

Moreover as $\M= B(x, 1),$ we have $1= |B(x,1)| \leq  (\frac 2\delta)^D |B(x,\delta)|$. Hence, 
\begin{equation}\label{inv}
 \hbox{ for all } x \in \M, \; \hbox{ for all }  0< \delta <1 \quad \frac 1{|B(x, \delta)|} \leq (\frac 2\delta)^D
\end{equation}

If $\M$ is connected one can prove additionally (see \cite{CKP}) that there exist $c>0, \; \beta >0,$ such that :
\begin{equation}\label{inv2}
 \hbox{ for all } x \in \M, \; \hbox{ for all }  0< \delta <1 \quad   c (\frac 1\delta)^\beta \leq \frac 1{|B(x, \delta)|} \leq (\frac 2\delta)^D
\end{equation}

\end{rem}
\subsection{Heat kernel} \label{sec-heat}

For this section, we follow standard expositions for heat kernel theory: for more details see
\cite{Ouhabaz}, \cite{Saloffbook}, \cite{Grigoryan}. 
We suppose that there exists  a self adjoint positive operator $L$ defined on a domain $D \subset \L^2$ dense in $ \L^2$.
Then $-L$ is the infinitesimal generator
of a self adjoint positive semigroup $e^{-tL}.$

We suppose in addition that $e^{-tL}$ is a Markov kernel operator i.e. there exists a  non negative kernel $ P_t(x,y)$ (\textit{'the heat kernel'}) such that :
\begin{align} &e^{-tL} f(x) = \int_{\M} P_t(x,y) f(y) d\mu(y)
\\
&P_t(x,y)=P_t(y,x), \\ &\int_{\M} P_t(x,y) d\mu(y) =1, \\
&P_{t+s}(x,y)= \int_{\M} P_t(x,u) P_s(u,y)du 
\end{align}

The following additional assumptions are central in our setting: there exist $C_1>0, c_1 >0, \alpha>0$, {such that }
$ \hbox{ for all } t \in ]0,1[, $

\begin{align} & 0\leq P_t(x,y) \leq \frac{C_1}{\sqrt{ |B(x, \sqrt t)| |B(y, \sqrt t)|}} e^{-\frac{c_1 \rho^2(x,y)}t} \label{bornesupgauss}
\\
  &|P_t(x,y)-P_t(x,y')|  \leq  \left[ \frac{\rho(y,y')}{\sqrt t}\right]^\alpha\frac{C_1}{\sqrt{ |B(x, \sqrt t)| |B(y, \sqrt t)|}} e^{-\frac{c_1 \rho^2(x,y)}t}
  \end{align}
 One can easily prove (see \cite{CKP})  that under the assumptions above, necessarily there exists  positive constants
$ C_1,\; C'_1$ such that
\begin{equation}\label{hker}
\frac{C'_1}{|B(x, \sqrt t)|}  \leq P_t(x,x) \leq \frac{C_1}{|B(x, \sqrt t)|}
\end{equation}

\begin{rem}
Under mild additional conditions on the space $\M$, see \cite{Ouhabaz}, Section 7.8, one can actually prove that  there exist $C_1, C_2>0, c_1, c_2 >0, \alpha>0$, {such that } $ \hbox{ for all } t \in ]0,1[, $
\[ \frac{C_2}{\sqrt{ |B(x, \sqrt t)| |B(y, \sqrt t)|}} 
e^{-\frac{c_2 \rho^2(x,y)}t}
\leq P_t(x,y) \leq \frac{C_1}{\sqrt{ |B(x, \sqrt t)| |B(y, \sqrt t)|}} 
e^{-\frac{c_1 \rho^2(x,y)}t}.
\]
On the last display, one can actually see that the heat kernel has a behaviour of square-exponential type. Furthermore, it is a positive definite kernel, see below.
\end{rem}

\subsection{Spectral decomposition}
The assumptions above have, as a consequence
(\cite{CKP}, Proposition 3.20) that the
 spectral decomposition of $L$ is discrete: there exists  a sequence $0=\lambda_0<\lambda_1 <\lambda_2< \ldots$ of
eigenvalues of $L$  associated with finite dimensional eigenspaces $\cH_{\lambda_k}$ such that :
$$ \L^2=  \oplus_{k} \cH_{\lambda_k}.$$ 
Necessarily $\cH_{\lambda_k} $ is a subset of $ C(\M)$ the space of continuous functions on 
$\M$.
More precisely, the projectors $P_{ \cH_{\lambda_k}}$ are kernel operators $P_k(x,y)$ with the following description: 
$$ P_k(x,y)= \sum_{1 \leq l \leq dim(\cH_{\lambda_k} ) } e_k^l(x)  e_k^l(y),$$
as soon as $\{e_k^l,\; 1 \leq l \leq dim(\cH_{\lambda_k} ) \} $ is an orthonormal basis of $\cH_{\lambda_k}.$
The Markov kernel $P_t$ writes:
\begin{equation} P_t(x,y)= \sum_{k} e^{-t \lambda_k} P_k(x,y)\label{markernel}
\end{equation}
Moreover, $e^{-tL}$ is a trace class operator. 
Using Mercer theorem, one can prove in addition that the convergence in the series is uniform.
\subsection{  Smooth functional calculus and `sampling-father-wavelets'} \label{smoothoperator}

More generally for any (very regular) $\Phi \in \D(\R)$ and $0<\delta \leq 1$, $\Phi(\delta \sqrt L)$ is a kernel operator described using  the spectral decomposion, via the following formula:
\begin{equation}
 \Phi(\delta \sqrt L)(x,y)= \sum_{k} \Phi(\delta \sqrt{\lambda_k}) P_k(x,y)\label{phi-delta}
\end{equation}
Our previous assumptions have the following consequences which will be important in the sequel, see \cite{CKP}:

\begin{description}
\item[{\it Localization:}] (\cite{CKP}, Section 3)
There exists a constant $C(\Phi)$ such that 
 \begin{equation}\label{FC}
    \hbox{ for all }\; 0<\delta \leq 1, \foral x,y\in \M,\; |\Phi(\delta \sqrt L)(x,y)| \leq   \frac 1{|B(x,\delta)| }   \frac{ C(\Phi)}{(1+ \frac{ \rho(x, y)}\delta )^{D+1}}.
    \end{equation}
%
From \eref{FC}  one can easily deduce the symmetrical bound
 $|\Phi(\delta \sqrt L)(x,y) |\leq \frac 1{\sqrt{|B(x,\delta)| |B(y,\delta)|    }}   \frac{ C(\Phi)}{(1+ \frac{ \rho(x, y)}\delta )^{D+1}}.$

\item[{\it Father wavelet:}] (\cite{CKP},  Lemmas 5.2 and 5.4)
There exist $ 0<C_0 <\infty,\quad 0<\gamma $
structural constants such that for any
$0<\delta \leq 1$ , for any  $\Lambda_{\gamma \delta}$ maximal $\gamma \delta-$net, there exists a family of functions :
$(D^\delta_\xi)_{\xi \in  \Lambda_{\gamma \delta}} $ such that 
\begin{equation}\label{wav0}
|D^\delta_\xi(x)| \leq  \frac 1{|B(x,\delta)|} \frac{ C_0}{(1+ \frac{ \rho(x, \xi)}\delta )^{d+1}},\;\quad \foral x\in \M
\end{equation}
and if we define the `low frequency' functions
\[\Sigma_t=\bigoplus_{\lam\le\sqrt{t}}\cH_\lam,\]
we have the following wavelet-type representation:
\begin{align}\label{wav}
&\foral \phi \in \Sigma_{1/\delta}, \quad \phi(x)=\sum_{\xi \in  \Lambda_{\gamma \delta}} \phi(\xi) |B(\xi, \delta)| D^\delta_\xi(x),
\\
\label{lemf} 
&\foral  (\alpha_\xi)_{\xi \in \Lambda_\delta },\; \| \sum_{\xi \in \Lambda_\delta }  \alpha_\xi |B(\xi,\delta)|  D^\delta_\xi (x) \|_\infty \lesssim   \sup_{\xi \in \Lambda_\delta }  |\alpha_\xi |
 \end{align}

We see on the formulae \eref{wav} and \eref{lemf} that the functions $|B(\xi, \delta)|D^\delta_\xi$ behave like father-wavelets, with coefficients directly obtained by sampling. We will see in Appendix B that these functions play an important role for instance to bound the entropy of various functional spaces.

\end{description}
 In the same spirit, we can also define an analogue of the mother wavelet. Notice that this construction will not appear explicitely in our Bayesian setting but will be used in the proof, see Section \ref{proof-entropy}. 
  Let us fix  $\Phi \in \D(\R), \; 0\leq \Phi, \;  1=\Phi(x), \; \hbox{for} \; |x|\leq 1/2, \; supp(\Phi) \subset [-1,1],$ and let us define also :
  $$ \Psi(x) =\Phi(\frac x2)- \Phi(x).$$
  So 
  $$ 0\leq \Psi (x)\leq 1, \; supp (\Psi) \subset \{ \frac 12 \leq |x|\leq 2\} ;  \quad
   \hbox{ for all } \delta >0 ,  \; \quad 1\equiv \Phi(\delta x) + \sum_{j\geq 0} \Psi (2^{-j}\delta x).$$
  
  So 
  $$ f =  \Phi(\delta \sqrt L)f + \sum_{j\geq 0} \Psi (2^{-j}\delta \sqrt L)f, \quad $$
  $$\Phi(\delta \sqrt L)f (x)=\int_M \Phi(\delta \sqrt L)(x,y) f(y)d\mu(y) ; \quad
   \Psi(\delta 2^{-j} \sqrt L)f (x)=\int_M  \Psi(\delta 2^{-j} \sqrt L)(x,y) f(y)d\mu(y)  $$
 $$ \Phi(\delta \sqrt L)f  \in \Sigma_{\frac 1\delta} ; \quad  \Psi (2^{-j}\delta \sqrt L)f  \in \Sigma_{\frac{2^{j+1}}\delta} \cap 
 [\Sigma_{\frac{2^{j+1}}\delta}]^{\perp}$$
  
%

 \subsection{ Examples.}

\paragraph{Torus case.} Let $\M=\bS^1$ be the torus equipped with the normalized  Lebesgue measure.
$\M$ is parameterised by $[-\pi, \pi]$ with identification of $\pi$ and $-\pi.$ The spectral decomposition of the Laplacian operator $\Delta$ gives rise to the classical Fourier basis, with 
$$\cH_0 =  span\{ 1\} ; \;  \ \cH_k= span\{e^{ikx}, e^{-ikx}\} = span\{ \sin kx, \cos kx\}$$
Hence,   $$dim (\cH_0)= 1; \; \hbox{ for all }  k >1, ~Ê dim (\cH_k) =2  \quad \hbox{and }\quad P_k(x,y) =2 \cos k(x-y).$$
$$ e^{t\Delta}(x,y) =1 + \sum_{k\geq 1} e^{-k^2t} 2\cos k(x-y) = \sqrt{\frac \pi t} \sum_{l \in \Z} e^{-\frac{(x-y-2l\pi)^2}{4t} }.$$
Clearly, $ \hbox{ for all } t >0, \quad  e^{t\Delta}(x,x) \geq 1.$  It holds 
$$ \hbox{ for all } 0<t<1, \quad  x,y \in [-\pi, \pi], \quad  
 C' \frac 1{\sqrt t} e^{-c'\frac{ \rho(x,y)^2}t}  \leq e^{t\Delta}(x,y)  \leq C \frac 1{\sqrt t} e^{-c\frac{ \rho(x,y)^2}t}.$$
Here we have, for any $x,y$ in $[-\pi, \pi]$, $$ \rho(x,y) =|x-y|\wedge (2\pi-|x-y|); \quad \hbox{ for all } 0<r\leq \pi, \quad   |B(x,r)| = \frac r\pi. $$

\paragraph{Jacobi case.} Let us now take $\M=[-1,1]$ equipped with the measure  $\omega(x) dx$ with
 $\omega(x) =(1-x)^\alpha (1+x)^\beta,~\alpha >- 1, ~\beta >-1.$ 
\\
If $\sigma(x)= (1-x)^2$, then  $\tau:=\frac{(\sigma \omega)'}{\omega} 
$ is a polynomial of degree $1$, we put :
$$-L(f)= D_J(f) = \frac{(\sigma \omega f')'}{\omega} = \sigma f'' + \tau f' $$
The operator  $L$ is a nonnegative symmetric  (in $\L_2 (\omega(x) dx)$) second order  differential operator (here and in the sequel, $u'$ denotes the derivative of $u$).
\\
Using Gram Schmidt orthonormalisation (again, in $\L_2 (\omega(x) dx)$) of $\{x^k,\, k\in \N\}$ we get a family of orthonormal polynomials $\{\pi_k ,\, k\in \N\}$ called Jacobi polynomials,
  which coincides with the spectral decomposition of $D_J.$ It holds
  $$D_J \pi_k = [k(k-1)\frac{\sigma"}2 +k \tau' ] \pi_k :=\lambda_k\pi_k 
  = - k(k +1+\alpha + \beta)) \pi_k
  $$
 Then, for any 
  $k \in \N$, $\;\cH_{\lambda_k}=span\{\pi_k\}, \;   ~dim(\cH_{\lambda_k}) =1$ and
  $$ ~P_k(x,y) = \pi_k(x) \pi_k(y);
  \quad \lambda_k =- k(k+\alpha + \beta +1).$$
  $$ e^{-t L(x,y)} = c_{\alpha, \beta}+ \sum_{k\geq 1} e^{-t k(k +1+\alpha + \beta))} \pi_k(x) \pi_k(y), \quad  c_{\alpha, \beta}^2 \int_M \omega(x) dx =1.$$
If $$ \rho(x,y) = | \arccos x - \arccos y| = \arccos (xy + \sqrt{1-x^2}  \sqrt{1-y^2}) $$
then
$$  C' \frac 1{  \sqrt{ |B(x, \sqrt t)|  |B(y, \sqrt t)|}} e^{-c' \frac{\rho^2(x,y) }t  }  \leq e^{-t L(x,y)}  \leq C \frac 1{  \sqrt{ |B(x, \sqrt t)|  |B(y, \sqrt t)|}} e^{-c \frac{\rho^2(x,y) }t  }$$
But $$ \hbox{ for all } x \in [-1,1], \quad   0<r  \leq \pi, \quad    |B(x,r)| \sim  r ((1-x) \vee r^2)^{\alpha +1/2} ((1+x) \vee r^2)^{\beta +1/2}.$$

\paragraph{Sphere case.} Let now $\M=\bS^{n-1} \subset \R^n$. 
 The geodesic distance on $\bS^{n-1}$ is given by 
$$\displaystyle{\rho(x,y) = \cos^{-1} (\langle x, y \rangle), ~~\langle x, y \rangle =\sum_{i=1}^n x_i y_i.}$$
There is a natural measure  $\sigma$ on $\bS^{n-1}$ which is \textit{ rotation invariant}.
%
 There is a natural Laplacian on $\bS^{n-1}, ~\Delta_{\bS^{n-1}}=\Delta,$ which is a negative self-adjoint operator with the following spectral decomposition.
 
 If $\cH_k$ is the restriction to  $\bS^{n-1}$ of polynomials of degree k
which are  \textit{homogeneous  } (i.e. $P(x) =\sum_{|\alpha| =k} a_\alpha x^\alpha,~~
 \alpha =(\alpha_1, \ldots, \alpha_n) , ~|\alpha| =\sum \alpha_i , ~\alpha_i\in \N$) and
\textit{harmonic  }
(i.e. $\Delta P= \sum_{i=1}^n \frac{\partial^2 P}{\partial x_i^2} =0,$)
we have,
$$ P \in \cH_{\lambda_k}(\bS^{n-1})  \implies  \ \ -\Delta_{\bS^{n-1}} P = k(k+ n-2) P:={\lambda_k}P$$
Moreover $  dim(\cH_k(\bS^{n-1})) = \frac{2k+n-2}{n-2} C^k_{n+k-3}= N_k(n).$
The space $\cH_k$ is called the space of spherical harmonics of order $k$,  (and with a slight abuse of notation  $\cH_k$ is our 'former'  $\cH_{\lambda_k}$). Moreover, if $(ÊY_{ki})_{1\leq i\leq N_k}$ is an  {orthonormal basis of} $\cH_k$, the projector writes
 $$P_k(x,y) = \sum_{1\leq i\leq  N_k} Y_{ki}(x) Y_{ki}(y).$$
 Actually :
  $$( x,y ) \in \bS^{n-1} \times \bS^{n-1}, \quad  P_k(x,y) = \frac 1{|\bS^{n-1}|} (1+ \frac k\nu) G^\nu_k( \langle x,y \rangle);$$
  $$ \quad \nu= \frac d2-1 ; \quad |\bS^{n-1}| =\frac{2\pi^{n/2}}{\Gamma(n/2)},$$
and  $G_k^\nu$ is the Gegenbauer polynomial of index $\nu$ and degree $k, $ defined for instance by its generating function :
$$ (1-2xt +t^2)^{-\nu} = \sum_k G_k^\nu(x) t^k.$$
Also, it holds
$$ |B(x,r)| =  |\bS^{n-2}| \int_0^r (\sin t)^{n-2} dt, $$ so at least for $ 0 \leq r \leq \pi/2$
$$  (\frac 2\pi)^{n-2}  \frac{ |\bS^{n-2}|}{n-1} r^{n-1} \leq |B(x,r)| \leq \frac{ |\bS^{n-2}|}{n-1} r^{n-1}$$
and clearly 
$$ \hbox{ for all }  0\leq r \leq \pi= diam(\bS^{n-1}), \quad c_1r^{n-1} \leq  |B(x,r)| \leq c_2r^{n-1}.$$
Now for all  $0\leq t \leq 1$, it holds
$$  \frac{C'}{ t^{(n-1)/2}} e^{-c' \frac{\rho(x,y)^2}t}  \leq e^{t \Delta}(x,y)= \sum_{k} e^{-t k(k+n-2)} P_k(x,y)  \leq \frac{C}{ t^{(n-1)/2}} e^{-c \frac{\rho(x,y)^2}t}.
$$

\paragraph{Ball case.}
Let $\M=\B^d$ be the unit ball of $\R^n.$  Let us consider the measure :
$$W(x) dx, \; W(x) =( 1-\| x\|^2)^{\mu-1/2}, \; \mu >0.$$
 Further define  the  operator
\begin{align*}Lf(x)&=
\Delta f(x)  -  {x}. \nabla({x}. \nabla f(x)) - (2\mu +d-1) 
{x}. \nabla f(x)\\
&=  \frac 1{W(x)}  div [(1-\| x\|^2) W(x)\nabla f(x)] + \frac 12\sum_{i\neq j} D^2_{i,j}f(x) , \quad \hbox{with}\quad
D_{i,j}f(x)= (x_j \partial_i -x_i \partial_j) f(x).
\end{align*}
One can verify that:
$$\int_{\B^d}  L(f) (x) f(x)W(x) dx
=
 - \int_{\B^d}  (1-\| x\|^2) |\nabla f |^2(x) W(x) dx - \frac 12\sum_{i\neq j} \int_{\B^d}  [D_{i,j}f ]^2(x)  W(x) dx.$$
Let $\Pi_k(\B^d)$ be the space of polynomials of degree  at most $k$ on the unit ball of $\R^d$
and define $ \V_k(\B^d) $ by
$$\Pi_k(\M)  = \V_k(\B^d) \bigoplus \Pi_{k-1}(\B^d) ;\hbox{ then}  $$
$$\L^2(\B^d) =   \bigoplus_{k=0}^\infty \V_k(\B^d).$$
The space $ \V_k(\B^d)$ is an eigenspace of $L.$ More precisely,
  $$f \in \V_k(\B^d)   \Longleftrightarrow L(f) = -k(k+d) f.$$
  The projector on $ \V_k$ is given by the following formula (see \cite{XUDUNKL}, \cite{PenchoXu}), with $\lambda = \mu + \frac{d-1}{2}$, 
$$  P_n(x,y)= c(d, \mu) \frac{n+\lambda}{\lambda}
\int_{-1}^1 G_n^\lambda \left(\langle x,y\rangle +
u\sqrt{1-|x|^2} \sqrt{1-|y|^2}\right) (1-u^2)^{\mu-1}du.$$
  $$
e^{-tL}(x, y) =\sum_{n\ge 0} e^{tn(n+\lambda)}P_n(x,y),
$$
The natural associated metric  is 
$$\rho_{\M}(x,y)= \arccos (\langle x,y \rangle + \sqrt{1-\| x\|^2} \sqrt{1-\| y\|^2} ).$$
Moreover, it holds
 $$| e^{tL} (x,y)| \leq C(c)\frac 1{ r^d (r+ \sqrt{1- \| x\|^2})^{\mu}(r+ \sqrt{1- \| y\|^2})^{\mu}}
e^{- c \frac{\rho_{\M}^2(x,y)}t}.$$
  But 
  $$ |B(x,r)| \sim r^d (r+ \sqrt{1- \| x\|^2})^{2\mu}
\sim r^d (r^2+ (1- \| x\|^2))^{\mu}  \sim r^d (r^2\vee (1- \| x\|^2))^{\mu}.$$
Clearly,
$$ |B(x,2r)| \leq 2^\delta |B(x,r)|, \quad \delta = d+ 2\mu.$$ 

\paragraph{Compact Riemannian manifold, without boundary.}
Let $\M$ be a compact Riemannian manifold of dimension $n$. Associated to the Riemannian structure we have a 
measure $dx$, a metric $\rho,$ and a Laplacian $\Delta.$
One has :
$$ \int_{\M} \Delta f(x) g(x) dx = -\int \nabla(f)(x). \nabla(g)(x) dx$$
So $-\Delta$ is a symmetric non negative operator.
Now the associated semigroup $e^{t\Delta}$ is a positive kernel operator verifying :
$$  C' \frac 1{  \sqrt{ |B(x, \sqrt t)|  |B(y, \sqrt t)|}} e^{-c' \frac{\rho^2(x,y) }t  }  \leq e^{t\Delta}(x,y) \leq C \frac 1{  \sqrt{ |B(x, \sqrt t)|  |B(y, \sqrt t)|}} e^{-c \frac{\rho^2(x,y) }t  }$$
The main property is that, see the appendix in Section \ref{appendix-manifold},
$$\exists  0< c_1 <c_2 <\infty , \quad \hbox{such that }, \quad \hbox{ for all } r \leq diam(M), \; c_1r^n \leq |B(x,r)| \leq c_2  r^n$$

\section{RKHS and heat kernel Gaussian process}\label{section_rkhs}
It is well known that for any $t>0$, $P_t(x,y)$ is  associated to a RKHS $\bH_t$
and there is a Gaussian  centered process $(W^t(x))_{x\in M}$ such that $E(W^t(x) W^t(y)) =P_t(x,y)$ for any $x,y$ in $\M$.
For instance,  $W^t(x)$ can be built in the following way :
$$ W^t(x) =\sum_i X_i \phi_i(x),$$ where $\phi_i(\cdot)$ is any orthonormal basis of $\bH_t,$
and  $\{X_i, \; i \in \N\}$ is a family of independent  Gaussian variables with mean 0 and variance 1. Also, $\bH_t$ is the isometric image of $\L^2$ by $P_t^{1/2}=P_{t/2}$. So the family $\{e^{-\lambda_kt/2}e^l_k,\;{k\in \N, \; 1 \leq l \leq dim(\cH_{\lambda_k}) }\}$ is a `natural' orthonormal basis of 
$\bH_t$. 
Hence,  we have the following description of $W_t$:
$$ W^t(x) = \sum_k \sum_{1\leq l \leq  \dim \cH_k} e^{-\lambda_kt/2}  X_k^l e^l_k(x),$$
where $X_k^l$ is a family of independent  Gaussian variables with variance 1.
\\
The RKHS $\bH_t$ has also the following description:
\[ \bH_t = \{h^t = \sum_k \sum_{1\leq l\leq dim(\cH_k)} a_k^l  e^{-\lambda_kt/2}e_k^l(x), \quad \sum_{k,l}|a_k^l|^2 <+\infty\},\]
equipped with the inner product 
\[ \langle \sum_k \sum_{1\leq l\leq dim(\cH_k)} a_k^l  e^{-\lambda_kt/2}e_k^l , 
 \sum_k \sum_{1\leq l\leq dim(\cH_k)} b_k^l  e^{-\lambda_kt/2}e_k^l \rangle_{\bH_t} 
 = \sum_k \sum_{1\leq l\leq dim(\cH_k)} a_k^l  b_k^l.\]
\\
Hence, if we denote by $\bH^1_t$ the unit ball of $\bH_t$:
$$ f\in \bH^1_t \Longleftrightarrow f=\sum_k \sum_{1\leq l\leq dim(\cH_k)} a_k^l  e^{-\lambda_kt/2}e_k^l(x), \quad \sum_{k,l}|a_k^l|^2 \leq 1.$$
\subsection{Entropy of the RKHS}
Let $(X, \rho)$ be a metric space. For $\epsilon >0$, we define,  as usual,  the covering number $N(\epsilon, X) $ as the smallest number of balls of radius
$\epsilon $ covering $X$. The entropy $H(\epsilon, X)$ is by definition $H(\epsilon, X)= \log_2 N(\epsilon, X) .$ \\

An important result of this section  is the link between the covering number $N(\epsilon, \M,\rho) $ of the space $\M$,     and 
$H(\epsilon, \bH^1_t, \L^p )$ for $p=2,\; \infty$ where $ \bH^1_t$ is the unit ball of the RKHS  defined in the subsection above.
More precisely we prove in appendix B, the following theorem:
\begin{theorem} \label{entropyconnection} Let us fix $\mu >0, \; a>0$. There exists $\epsilon_0>0$ such that for $\epsilon, t$ with $ \epsilon^\mu \leq at$
and  $ 0<\epsilon\leq \epsilon_0,$ 
$$ H(\epsilon, \bH^1_t, \L^2) \sim   H(\epsilon, \bH^1_t, \L^\infty) \sim    N(\delta(t,\epsilon), \M)  . \log \frac 1\epsilon
\quad \hbox{where}\quad \frac 1{\delta(t,\epsilon)}:=
  \sqrt{\frac 1t  \log (\frac 1\epsilon)}.$$
\end{theorem}
\begin{rem} \label{rementropy}
Theorem \ref{entropyconnection} gives the precise behaviour up to constants, from above and below, of the entropy of the RKHS unit ball $\bH_t^1$. The constants involved depend only on $\M, a, \mu$. The mild restriction on the range of $t$ arises from technical reasons in the proof of the upper-bound. As an examination of the proof reveals, this restriction is not needed in the proof of the lower bound.
\end{rem}

\subsection{Uniform  polynomial control for the measures of the balls }
In Appendix B, the general case is considered, but  for sake of simplicity,  in the sequel, we will concentrate on the following case where the entropy of $\M$ has an exact polynomial control.
\\
As a matter of fact, in number of examples, for instance for compact Riemannian manifolds without boundary, see Appendix C, the bounds in \eref{inv2} are the same and we have the following uniform polynomial control (Ahlfors condition, see for instance \cite{juha01}). 
 There exist $  c_1>, c_2>0,  d>0 $ such that 
\begin{equation}\label{poly-case} \hbox{ for all } x \in \M, \;  \hbox{ for all } \;  0 <r  \leq 1, \; c_1 r^d \leq | B(x,r)| \leq c_2 r^d
\end{equation} 
Necessarily, $ d \leq D$ since using \eref{inv}, we have $|B(x,r)| \geq (r/2)^D$.
\\
In this case, Theorem \ref{entropyconnection} takes the following form.

\begin{proposition}
In the case \eref{poly-case}, we have
\begin{align}
\label{poly-entropie1} \frac{1}{c_2}  (\frac 1\epsilon)^d  &\leq N(\epsilon, \M) \leq  card(\Lambda_{\epsilon})  \leq  \frac{2^d}{c_1} (\frac 1\epsilon)^d\\
 H(\epsilon, \bH^1_t, \L^2) &\sim   H(\epsilon, \bH^1_t, \L^\infty) \sim  (\frac 1{ \delta(t,\epsilon)})^d \log \frac 1\epsilon.\label{poly-entropie2}
\end{align}
For all $ 0<\epsilon\leq \epsilon_0$; if we suppose that for some $\mu >0, \; a>0, \;   \epsilon^\mu \leq at.$

\end{proposition}
\noindent 
Indeed, let  $(B(x_i,\epsilon))_{i\in I}$ be a minimal covering of $\M$; we have
$$ 1= |\M| \leq \sum_{i\in I}  |B(x_i, \epsilon)|  \leq N(\epsilon, \M) c_2 \epsilon^d.$$
Now if  $\Lambda_{\epsilon}$ is any maximal $ \epsilon-$net, we have:
$$ 1= |\M| \geq \sum_{\xi \in \Lambda_{\epsilon}} |B(\xi, \epsilon/2) | \geq card(\Lambda_{\epsilon}) c_1 (\epsilon/2)^d.$$

\section{Geometrical Prior: concentration function} \label{sec-prior}
For the Bayesian results in the next four sections, it is assumed that the compact metric space $(\M,\rho)$ is as in Section \ref{sec-metspace}, that the polynomial estimate \eqref{poly-case} for volume of balls holds and
that there exists a heat-kernel operator with the properties listed in Section \ref{sec-heat}.

\subsection{Prior, definition} We consider a prior on functions constructed hierarchically as follows. First draw a positive random variable $T$ according to a density $g$ on $(0,1]$. Suppose that \eqref{poly-case} holds and assume there exists a real $a>1$ and positive constants $c_1,c_2,q$ such that, with $d$   defined in \eref{poly-case},
\begin{equation} \label{priorT}
 c_1 t^{-a} e^{-t^{-d/2}\log^q(1/t)}    \le g(t) \le  c_2 t^{-a} e^{-t^{-d/2}\log^q(1/t)}, \quad t\in (0,1].
\end{equation}
We show below that the choice $q=1+d/2$ leads to sharp rates. The choice of this particular form for the prior on $t$ is related to the form taken by the entropy of $\bH_t^1$. For more discussion on this, see Section \ref{sec-lb}. Also, note that we do not consider large values of $t$. This correspond to the fact that the trajectories of the process $W^t$ are already very smooth (in fact, infinitely differentiable almost surely). So, to capture rates of convergence for typical smoothness levels such as Sobolev or Besov indexes, we only need to make the paths `rougher', which corresponds to $t$ small.

Given $T=t$, generate a collection of independent standard normal variables 
$\{X_k^l \}$ with indexes ${k\ge 0, \; 1 \leq l \leq dim(\cH_{\lambda_k}) }$ and set, for $x$ in $\M$, 
\begin{equation} \label{prior}
  W^t(x) = \sum_k \sum_{1\leq l \leq  \dim \cH_k} e^{-\lambda_kt/2}  X_k^l e^l_k(x).  
\end{equation}
In the following the notation $W^t$ refers to the Gaussian prior \eqref{prior} for a fixed value of $t$ in $(0,1]$. One can check that $W^t$ defines a Gaussian variable in various separable Banach spaces $\mb$, see \cite{rkhs} for definitions. More precisely, we focus on the cases $\mb=(\C^0(\M),\|\cdot\|_{\infty})$ and $ \mb=(L^2(\M,\mu),\|\cdot\|_2)$.
To do so, apply Theorem 4.2 in \cite{rkhs}, where almost sure convergence of the series \eqref{prior} in $\mb$ follows from the properties of the Markov kernel \eqref{markernel}.

 The full (non-Gaussian) prior we consider is $W^T$, where $T$ is random with density $g$ given by \eqref{priorT}.
Hence, this construction leads to a prior $\Pi_w$, which is the probability measure induced by
\begin{equation} \label{priorw}
 W^T(x) =  \sum_k \sum_{1\leq l \leq  \dim \cH_k} e^{-\lambda_kT/2}  X_k^l e^l_k(x). \end{equation} 
\subsection{Approximation and small ball probabilities}
The so-called concentration function of a Gaussian process defined below turns out to be fundamental to prove sharp concentration of the posterior measure. For this reason  we focus now on the detailed study of this function for the geometrical prior.


In this paper, the notation $(\mb,\|\cdot\|_\mb)$ is used for anyone of the two spaces
\[(\mb,\|\cdot\|_\mb) = (\C^0(\M),\|\cdot\|_\infty) \quad\text{or}\quad
 (\mb,\|\cdot\|_\mb) = (\L^2,\|\cdot\|_{2}). \]
Any property stated below with a $\|\cdot\|_{\mb}$-norm holds for both spaces.

\paragraph{Concentration function.}  
Consider the Gaussian process \eqref{prior} $W^t$, for a fixed $t\in (0,1]$. Its concentration function within $\mb$ is defined, for any function $w_0$ in $\mb$, as the sum of two terms
\begin{align*} 
 \varphi_{w_0}^t(\eps) & =
 \inf_{h^t\in \bH_t,\ \|w_0-h_t\|_{\mb} < \eps} \frac{\| h_t\|_{\bH_t}^2}{2}
- \log \P( \|W^t\|_\mb <\eps) \\
 & := \qquad A_{w_0}^t(\eps)  \qquad +  \qquad S^{t}(\eps). 
\end{align*}
Notice that the approximation term quantifies how well $w_0$ is approximable by elements of the RKHS $\bH_t$ of the prior while keeping the `complexity' of those elements, quantified in terms of RKHS-norm, as small as possible. The term $A_{w_0}^t(\eps)$ is finite for all $\eps>0$ if and only if $w_0$ lies in the closure in $\mb$ of $\bH_t$ (which can be checked to coincide with the support of the Gaussian prior, see \cite{rkhs}, Lemma 5.1.) It turns out that for the prior $W^t$ this closure is $\mb$ itself, as quite directly follows from the approximation 
results below.

In order to have a precise calibration of $A_{w_0}^t(\eps)$, we will assume regularity conditions on the function $w_0$, which in turn will yield the rate of concentration.
Namely we shall assume that $w_0$ belongs to a regularity class $\F_s(\M)$, $s>0$ taken equal to 
a Besov space
\[ \F_s(\M) = B_{\infty,\infty}^s(\M)\ \  \text{if}\quad \mb
= \C^0(\M) \qquad  \text{(resp.}\  \ B_{2,\infty}^s(\M)\ \  \text{if}\ \ \mb
= \L^2 ).\]
The problem of the regularity assumption in a context like here is not a simple one. We took here a natural generalization of the definition of usual spaces on the real line, by means of approximation property. For more details we refer to
 Appendix A. It holds if $\Phi$ is a Littlewood-Paley function (i.e. verifies the 
conditions \eref{cond_phi} of Appendix A) and $w_0 \in \F_s(\M)$,

{ \[ \|\Phi(\delta \sqrt{L})w_0 - w_0 \|_{\mb} \le C \delta^s =: \eps.\]}

\paragraph{Approximation term $A_{w_0}^t(\eps)$-regularity assumption on $\M$.  } 

For any $w_0$ in the Banach space $\mb$, consider the sequence of approximations, for $\delta\to 0$ and $L$ defined above, using \eref{phi-delta}, with a function $\Phi$ having support in $[0,1]$,
\[ \Phi(\delta \sqrt{L}) w_0 = \sum_{\la_k \le \delta^{-2}} \Phi(\delta\sqrt{\la_k})P_{\cH_{\la_k}} w_0, \]
where $P_{\cH_{\la_k}}$ is the projector onto the linear space spanned by the $k$th-eigenspace $\cH_{\lambda_k}$
defined above. For any $\delta>0$, the sum in the last display is finite thus 
$\Phi(\delta \sqrt{L}) w_0$ belongs to $\bH_t$.

On the other hand, making use of the previous choice $\delta^s =: \eps$,
{ \begin{align*}
 \| \Phi(\delta \sqrt{L}) w_0 \|_{\bH_t}^2 & 
 = \sum_{\la_k \le \delta^{-2}} |\Phi(\delta\sqrt{\la_k})|^2 
 e^{\la_k t} \| P_{\cH_{\la_k}} w_0 \|_2^2 \\
& \le C \sum_{\la_k \le \delta^{-2} }  e^{\la_k t}  \| P_{\cH_{\la_k}} w_0 \|_2^2 
\\
&  \le C e^{t\eps^{-2/s} }\|  w_0 \|_2^2. 
\end{align*}}
Note that $\|  w_0 \|_2\le 1$ if we suppose that $w_0$ is in the unit ball of  $\F_s(\M)$ (since necessarily $\|w_0\|_\mb$ is bounded by 1 and, for the case of the infinity norm, since $\M$ is compact with  $\mu$-measure $1$).
Hence, we proved
\begin{equation}\label{approxterm-upper}
A_{w_0}^t(\eps)\le C e^{t\eps^{-2/s} },\quad \hbox{ if } w_0 \in \F_s(\M).
\end{equation}
Note that this is precisely the place where the regularity of the function plays a role.
\paragraph{Small ball probability $S^t(\eps)$.}

Let us show in successive steps that the following upper-bound on the small ball probability of the Gaussian process $W^t$ viewed as a random element in $\mb$ holds.
\begin{proposition}
Fix $A>0$. There exists a universal constant $\eps_0>0$, and constants $C_0,C_1>0$ which depend on $d,A, \mb$ only, such that, for any $\eps\le\eps_0$ and any $t\in[C_1\eps^{A},1]$,
 \begin{equation}\label{smallball}
-\log \P\left( \| W^t \|_{\mb} \le \eps \right) \le C_0 
\left(\frac{1}{\sqrt{t}}\right)^d \left(\log \frac{1}{\eps}  \right)^{1+\frac{d}{2}}.\\ 
\end{equation}
\end{proposition}


{\em The steps follow the method proposed by \cite{vvvz09}. }
The starting point is a bound on the entropy of the unit ball $\bH_t^1$ of $\bH_t$ with respect to the sup-norm, which is a direct consequence of \eref{poly-entropie2} and is summarized by the following:


There exists a universal constant $\eps_1>0$, and constants $C_2,C_3>0$ which depend on $d,A$ only, such that, for any $\eps\le\eps_1$ and any $t\in[C_2\eps^{A},1]$,

\begin{equation} \label{entropy}
 \log N(\eps,\bH_t^1,\|\cdot\|_\mb) \le C_3 \left(\frac{1}{\sqrt{t}}\right)^d \left(\log \frac{1}{\eps}  \right)^{1+\frac{d}{2}}.\\
\end{equation}

{\bf Step 1, crude bound.} $\ $
Let $u_t$ be the mapping canonically associated to $W^t$ considered in \cite{LiLinde99} and, as in this article, set
\begin{align*}
e_n(u_t) & :=\inf\ \{\eta>0,\ N(\eta,\bH_t^1,\|\cdot\|_\mb) \le 2^{n-1} \} \\
 & \le \inf \{0<\eta<t,\ \log N(\eta,\bH_t^1,\|\cdot\|_\mb) \le (n-1)\log 2 \}. 
\end{align*}
By definition, the previous quantity is smaller than the solution of the following equation in $\eta$, where we use the bound \eqref{entropy},
\[ C t^{-\frac{d}{2}}\log^{1+\frac{d}{2}} \frac{1}{\eta} = n  \] 
that is $\eta = \exp\{-Cn^{\frac{2}{2+d}} t^{\frac{d}{2+d}} \}$. Thus
\[ e_n(u_t) \le \exp\{-Cn^{\frac{2}{2+d}} t^{\frac{d}{2+d}} \},\quad n\ge 1.\]
The first equation of \cite{Tomczak87}, page 300 can be written
\[ \sup_{k\le n} k^\al e_k(u^*_t) \le 32 \sup_{k\le n} k^\al e_k(u_t).\]
We have, for any $k\ge 1$ and any $m\ge 1$,
\begin{align*}
k^{2m} e_k(u_t) & \le k^{2m} \exp\{-Ck^{\frac{2}{2+d}} t^{\frac{d}{2+d}}\} \\
 & \le  t^{-md} (k^{2}t^{d})^m \exp\{-C (k^2 t^d)^{\frac{1}{2+d}} \} \\
 & \le t^{-md} V_m(k^{2}t^{d}),
\end{align*}
where $V_m:x\to x^m e^{-Cx^{\frac{1}{2+d}}}$ is uniformly bounded on $(0,+\infty)$ by a finite constant $c_m$ (we omit the dependence in $d$ in the notation). It follows that for any $n\ge 1$,
\begin{align*}
n^{2m} e_n(u^*_t) & \le \sup_{k\le n} k^{2m} e_k(u^*_t) \\
 & \le 32 \sup_{k\le n} k^{2m} e_k(u_t) \\
 & \le 32 c_m t^{-md}.
\end{align*}
We have obtained $e_k(u_t) \le 32 c_m t^{-md} k^{-2m}  $ for any $k\ge 1$.
Lemma 2.1 in \cite{LiLinde99}, itself cited from \cite{Pisier89}, can be written as follows. If $\ell_n(u_t)$ denotes the $n$-th approximation number of $u_t$ as defined in \cite{LiLinde99} p. 1562,
\[ \ell_n(u_t) \le c_1 \sum_{k\ge c_2 n} e_k(u^*_t) k^{-1/2}(1+\log{k}). \]
From the bound on $e_k(u^*_t) $ above one deduces, for some constant $c_m'$ depending only on $m$, for any $n\ge 1$,
\[ \ell_n(u_t) \le c_m' t^{-d} n^{1-2m}.\]

Consider the definitions, for any $\eps>0$ and $t>0$,
\[ n_t(\eps):=\max\{n:\ 4\ell_n(u_t)\ge \eps\},\qquad 
   \sigma(W^t) = \E\left[ \|W^t\| \right]^{1/2}. \] 
 A sufficient condition for $n_t(\eps)$ to exist is $4\sigma(W^t) \ge \eps$, since
$\ell_n(u_t)\le \ell_1(u_t)=\sigma(W^t)$. So, provided $\eps\le 4\sigma(W^t)$,  the 
bound on $\ell_n$ implies 
$n_t(\eps)\le C_m (\eps^{-1} t^{-d})^{1/(2m-1)}$.

The following result makes Proposition 2.3 in  \cite{LiLinde99} precise with respect to constants involving the process under consideration. This is important in our context since we consider a collection of processes $\{W_t\}$ indexed 
by $t$ and need to keep track of the dependence in $t$.

\begin{proposition} \label{prop23}
Let $X$ be centered Gaussian in a real separable Banach space $(E,\|\cdot\|)$. Define $n(\eps)$ and $\sigma(X)$ as above. Then for a universal constant $C_4>0$, any $\eps\le 1 \wedge (4\sigma(X))$,
\[ -\log\P\left[  \|X\|<\eps \right] \le 
C_4 n(\eps)\log \left[  \frac{6n(\eps)( \sigma(X) \vee 1)}{\eps} \right].
\]
\end{proposition}

\vspace{.5cm}

Explicit upper and lower bounds for $\sigma(W^t)$ are given in Appendix B, see 
\eref{w-infty}-\eqref{w-l2}.  In the `polynomial case', see \eref{poly-case}, these bounds imply, uniformly in the interval of $t$'s considered, that $1 \lesssim \sigma(W^t)\lesssim \eps^{-B}$ for some $B>0$,

Combining this fact with Proposition \ref{prop23} and the previous bound on $n_t$,  we obtain that for some positive constants $C_7,\eps_3, \zeta$, for any $\eps\le \eps_3$ 
and $t\in [C_2\eps^{A},1]$
\begin{equation} \label{crudebound}
S^t(\eps)=-\log \P\left( \| |W^t \|_\mb \le \eps \right) \le C_7 \eps^{-\zeta}.
\end{equation}

{\bf Step 2, general link between entropy and small ball.} $\ $
According to Lemma 1 in \cite{KuelbsLi93}, we have, 
if $G$ is the distribution function of the standard Gaussian distribution (see their formula (3.19), or (3.2)),
\[ S^t(2\eps) + \log G(\la+G^{-1}(e^{-S^t(\eps)}))
\le \log N\left(\frac{\eps}{\la},\bH_t^1,\|\cdot\|_\mb \right).\]
Lemma 4.10 in \cite{vvvz09} implies, for every $x>0$, 
\[ G(\sqrt{2x}+G^{-1}(e^{-x}))\ge 1/2.\] 
Take $\la=\sqrt{2 S^t(\eps)}$ in the previous display. Then for values of $t,\eps$ such that \eqref{entropy} holds,
\[ S^t(2\eps) + \log\frac{1}{2} \le 
C  \left(\frac{1}{\sqrt{t}}\right)^d 
\left(\log \frac{S^t(\eps)}{\eps}  \right)^{1+\frac{d}{2}}.\]
Finally combine this with \eqref{crudebound} to obtain the desired Equation \eqref{smallball} that is
\begin{equation*}
S^t(\eps) \le C 
\left(\frac{1}{\sqrt{t}}\right)^d \left(\log \frac{1}{\eps}  \right)^{1+\frac{d}{2}}.
\end{equation*}
under the conditions $\eps\le \eps_3$ and 
$C_2\eps^{A}\le t\le 1$.

\section{General conditions for posterior rates}\label{section_gc}

A general theory to obtain convergence rates for posterior distributions relative to some distances is presented in \cite{ggv00} and \cite{gvni}.
 The object of interest is a function $f_0$ (e.g. a regression function, a density function etc.). In some cases, for instance 
 density estimation with  Gaussian priors, one cannot directly put the prior on the density itself (a Gaussian prior does not lead to positive paths). This is why we will parametrize the considered statistical problem with the help of a function $w_0$ in some separable Banach space $(\mb,\|\cdot\|_\mb)$ of functions defined over $(\M,\rho)$.  In some cases (e.g. regression) $w_0$ and $f_0$ coincide, in others not (e.g. density estimation), see examples below. As before, $\mb$ is either $\C^0(\M)$  or $\L^2$. 
 \\
In this Section we check that there exist Borel measurable subsets $B_n$ in $(\mb,\|\cdot\|_\mb)$ such that, for some vanishing sequences $\eps_n$ and $\bar{\eps}_n$, some $C>0$ and $n$ large enough,
\begin{align}
\P(\|W^T-w_0\|_\mb \le \eps_n) & \ge e^{-Cn\eps_n^2}   \label{priormass}\\
\P(W^T\notin B_n) & \le e^{-(C+4)n\eps_n^2}         \label{sieve} \\
\log N(\enb, B_n,\|\cdot\|_\mb) & \le n \enb^2  \label{entropie}
\end{align}
In Section \ref{section_rates}, we show how this quite directly implies posterior concentration results.
In \cite{vvvz09}, the authors also follow this approach. One advantage of the prior considered here is that, contrary to \cite{vvvz09}, the RKHS unit balls are precisely nested as the time parameter $t$ varies, see \eqref{nested}. This leads to slightly simplified proofs.

\paragraph{Prior mass.}  For any fixed function $w_0$ in $\mb$ and any $\eps>0$, by conditioning on the value taken by the random variable $T$,
\begin{align*}
\P(\|W^T-w_0\|_\mb<2\eps) & = \int_0^1 \P(\|W^t-w_0\|_\mb<2\eps) g(t)dt.
\end{align*}
The following inequality links mass of Banach-space balls for Gaussian priors with their concentration function in $\mb$, see \cite{rkhs}, Lemma 5.3,
\[ e^{-\varphi^t_{w_0}(\eps)} \le 
 \P(\|W^t-w_0\|_\mb<2\eps) \le e^{-\varphi^t_{w_0}(2\eps)} ,\]
for any $w_0$ in the support of $W^t$. We have seen above that any $f_0$ in $\F^s(\M)$ belongs to the support of the prior. It is not hard to adapt the argument to check that in fact any $f_0$ in $\mb$ can be approximated in $\mb$  by a sequence of elements in the RKHS ${\bH}_t$ and thus belongs to the support in $\mb$ of the prior by Lemma 5.1 in \cite{rkhs}. Then 
\begin{align*}
\P(\|W^T-w_0\|_\mb<2\eps) & \ge \int_0^1 e^{-\varphi^t_{w_0}(\eps)} g(t)dt\\
& \ge \int_{t_\eps^*}^{2t_\eps^*} e^{-\varphi^t_{w_0}(\eps)} g(t)dt,
\end{align*}
for some $t_\eps^*$ to be chosen.

The concentration function is bounded from above, under the conditions  $\eps\le \eps_3 $ and  $t$ in $[C_2\eps^A,1]$,  by
\[ \varphi^{t}_{w_0}(\eps)
\le C \left[ e^{\eps^{-2/s} t} +  
\left(\frac{1}{\sqrt{t}}\right)^d \left(\log \frac{1}{\eps}  \right)^{1+\frac{d}{2}}\right] \]

Set $t_\eps^*=\delta \eps^{\frac{2}{s}}
\log\frac{1}{\eps}$ with $\delta$ small enough to be chosen. This is compatible with the above conditions provided $A>2/s$.
Then for $\eps$ small enough and any $t\in[t_\eps^*,2t_\eps^*]$, 
\[ \varphi^{t}_{w_0}(\eps)
\le C \left[ \eps^{-2\delta} +  
\delta^{-d}\eps^{-\frac{d}{s}}\left(\log\frac{1}{\eps}  \right)\right] . \]
Set $\delta=d/(4s)$.
One obtains, for any $t\in[t_\eps^*,2t_\eps^*]$, 
\[\varphi^{t}_{w_0}(\eps)
\le C_d \eps^{-\frac{d}{s}}\left(\log\frac{1}{\eps}  \right). \]
Inserting this estimate in the previous bound on the prior mass, one gets, together with \eqref{priorT}, for $\eps$ small enough and $q\le 1+d/2$,
\begin{align}
\P(\|W^T-w_0\|_\mb<2\eps) 
& \ge  t_\eps^* e^{-C\eps^{-\frac{d}{s}}\left(\log\frac{1}{\eps} \right) } 
\left[ \inf_{t\in[t_\eps^*,2t_\eps^*]} g(t) \right] \nonumber\\
& \ge C{t_\eps^*}^{1-a} e^{- {t_\eps^*}^{-\frac{d}{2}}\left(\log\frac{1}{t_\eps^*}\right)^q  - C\eps^{-\frac{d}{s}}
\left(\log\frac{1}{\eps}\right) }\nonumber\\
& \ge C {\eps}^{2(1-a)/s} (\log\frac{1}{\eps})^{1-a}
  e^{- C \eps^{-\frac{d}{s}}
\left(\log\frac{1}{\eps}\right)^{q-\frac{d}{2}} - C\eps^{-\frac{d}{s}}
\left(\log\frac{1}{\eps}\right) }\nonumber\\
& \ge C e^{ - C'\eps^{-\frac{d}{s}}\left(\log\frac{1}{\eps}\right) }. \label{pmcond}
\end{align}
Condition \eqref{priormass} is satisfied for the choice 
\begin{equation} \label{rate}
\eps_n \sim \left( \frac{n}{\log n}\right)^{-\frac{s}{2s+d}}.
\end{equation}

\paragraph{Sieve.}  The idea is to build sieves using Borell's inequality. 
Recall here that $\bH_r^1$ is the unit ball of the RKHS of the centered Gaussian
 process $W^r$, viewed as a process on the Banach space $\mb$. The notation $\mb_1$ (as well as $\bH_r^{1}$) stands for the unit ball of the associated space. 

First, notice that from the explicit form of the RKHS of $W^t$, we have
\begin{equation}\label{nested}
 \text{If }\ t_2\ge t_1,\qquad{then}\qquad \bH_{t_1}^{1}\subset\bH_{t_2}^{1}.
\end{equation}


Let us set for $M=M_n$, $\eps=\eps_n$ and $r>0$ to be chosen later, 
\[B_n = M\bH_r^{1} + \eps \mb_1, \]
Consider the case $t\ge r$, then using \eqref{nested}
\begin{align}
\P(W^t\notin B_n) 
& = \P(W^t \notin M\bH_r^{1} + \eps \mb_1) \nonumber \\
& \le \P(W^t \notin M\bH_t^{1} + \eps \mb_1)\nonumber \\
& \le 1- G(G^{-1}(e^{-S^t(\eps)}) + M).\label{eqborell}
\end{align}
where the last line follows from Borell's inequality.




{\bf Choices of $\eps$, $r$ and $M$.} $\ $ Let us set  $\eps=\eps_n$ given
 by \eqref{rate} and
\begin{equation}\label{choixrm}
r^{-\frac{d}{2}} \sim \frac{n\eps_n^2}{(\log n)^{1+\frac{d}{2}}}
 \qquad\qquad\text{and}\qquad\qquad M^2 \sim n\eps_n^2. 
\end{equation}
 First, one checks that $r$ belongs to $[C_2\eps^A,1]$. This is clear from the definition since we have assumed $A>2/s$.  
Then any $t\in[r,1]$ also belongs to $[C_2\eps^A,1]$ so we can use the entropy bound  and write
\[ S^t(\eps) \le C t^{-\frac{d}{2}}\left( \log \frac{1}{\eps_n} \right)^{1+\frac{d}{2}}
\le C r^{-\frac{d}{2}}\left( \log \frac{1}{\eps_n} \right)^{1+\frac{d}{2}}=:S^*_n. \]
Now the bounds $-\sqrt{2\log(1/u)}\le G^{-1}(u)\le -\frac{1}{2}\sqrt{\log(1/u)}$
valid for $u\in(0,1/4)$ imply that
\[1- G(G^{-1}(e^{-S^t(\eps)}) + M) \le 
1- G(G^{-1}(e^{-S^*_n}) + M) 
 \le e^{-M^2/8},\]
as soon as $M\ge 4\sqrt{S^*_n}$ and $e^{-S^*_n}<1/4$.

 To check $e^{-S^*_n}<1/4$ note that $S_n^*\ge S^r(\eps)$ which can be further bounded from below using Equation (3.1) in \cite{KuelbsLi93} which leads to, for any $\eps,\lambda>0$,
\[S^r(\eps) \ge H(2\eps,\la \bH_r^1) - \frac{\la^2}{2} \ge Cr^{-\frac{d}{2}}(\log\frac{\lambda}{\eps})^{1+\frac{d}{2}} - \frac{\lambda^2}{2}.\]
Here we have used the bound from below of the entropy see \eref{poly-entropie2}. Then take $\la=1$ to obtain $S^*_n(\eps)\ge \log(4)$ for $\eps$ small enough. 

The first inequality $M\ge 4\sqrt{S^*_n}$ is satisfied if
\[ M^2 \ge 16 r^{-\frac{d}{2}}\left( \log \frac{1}{\eps_n} \right)^{1+\frac{d}{2}},\]
and this holds for the choices of $r$ and $M$ given by \eqref{choixrm}.
Hence for large enough $n$,
\begin{align*}
\P(W^t\notin B_n) & \le e^{-M^2/8} \\
& \le e^{-Cn\eps_n^2}.
\end{align*}
Then we can write, if $q\ge 1+d/2$,
\begin{align*}
\P(W^T\notin B_n) 
& = \int_0^{1} \P(W^t\notin B_n)  g(t)dt \\
& \le \P(T< r) + \int_r^{1} \P(W^t\notin B_n)  g(t)dt \\
& \le Cr^{-c}e^{-C'r^{-d/2}\log^q(\frac{1}{r})} + e^{-M^2/8}\\
&  \le e^{-Cn\eps_n^2}.
\end{align*}

\paragraph{Entropy.} 
It is enough to bound from above
\begin{align*}
\log N(2\eps_n, M\bH_r^1+\eps_n\mb_1,\|\cdot\|_\mb) 
& \le \log N(\eps_n, M\bH_r^1,\|\cdot\|_\mb) \\
& \le r^{-d/2} \left(\log \frac{M}{\eps_n} \right)^{1+\frac{d}{2}}\\
& \le C n\eps_n^2, 
\end{align*}
where we have used \eqref{entropy} to obtain the one but last inequality.

\section{Posterior rate, main results } \label{section_rates}

In the next paragraphs, we recall the definition of the Bayes posterior measure, in a dominated setting where the posterior is given by Bayes' formula, and state a general rate-Theorem. We then prove the announced results in the three considered statistical settings. We study the convergence of the posterior measure in a frequentist sense in that we suppose that there exists a `true' parameter, here an unknown function, denoted $f_0$. That is, we consider convergence under the law of the data under $f_0$, and denote the corresponding distribution $P_{f_0}^{(n)}$. The expectation under this distribution is denoted $\E_{f_0}$.

For any densities $p, q$ with respect to a measure $\mu$, denote 
\[ K(p,q) = \int p\log{\frac{ p}{q}}d\mu,\qquad
V(p,q) = \int p\log^2{\frac{ p}{q}}d\mu.\]

\subsection{Bayesian framework and general result}
Let $\F$ be a metric space equipped with a $\sigma$-field $\cT$. \\
{\em Data.} Consider a sequence of statistical experiments 
$(\cX_n,\A_n,\{P_{f}^{(n)}\}_{f\in\F})$ indexed by the space $\F$. We assume that there exists a common ($\sigma$-finite) dominating measure $\mu^{(n)}$ to all probability measures $\{P_{f}^{(n)}\}_{f\in\F}$, that is
\[ dP_{f}^{(n)}(x^{(n)}) = p_f(x^{(n)}) d\mu^{(n)}(x^{(n)}). \]
We also assume that the map
$ (x^{(n)},f) \to  p_f(x^{(n)}) $ is jointly measurable relative to $\A_n \otimes \cT$. 
{\em Prior.} We equip the space $(\F,\cT)$  of a probability measure $\Pi$ that is called prior. Then the space $\cX_n\times \F$ can be naturally equipped of the $\sigma$-field $\A_n\otimes\cT$ and of the probability measure
\[ P(A_n\times T) = \int\int_{A_n\times T} p_f^{(n)}(x^{(n)})d\mu^{(n)}(x^{(n)})d\Pi(f).\]
The marginal in $f$ of this measure is the prior $\Pi$. The law $X\given f$ is $P_f^{(n)}$.  \\
{\em Bayes formula.} 
Under the preceding framework, the conditional distribution of $f$ given the data $X^{(n)}$ is absolutely continuous with respect to $\Pi$ and is given by, for any measurable set $T\in \cT$,
\[ \Pi(T\given X^{(n)}) = \frac{ \int_T  p_f(X^{(n)}) d\Pi(f) }{\int  p_f(X^{(n)}) d\Pi(f) }. \]

Let $\eps_n\to0$ be a rate of convergence such that $n\eps_n^2\to+\infty$.
Define a Kullback-Leibler neighborhood of the element $f_0$ in $\F$ by
\[ B_{KL}(f_0,\eps_n)=\
\Big\{\ K(p_{f_0}^{(n)},p_f^{(n)}) \le n\eps_n^2,\ \  V(p_{f_0}^{(n)},p_f^{(n)})  \le n\eps_n^2  
\ \Big\}   \]

Next we state a general result which gives sufficient conditions for the convergence of the posterior measure. It is a slight variation on Theorem 1 in \cite{gvni}.
A first key ingredient is the existence of a distance $d_n$ enabling testing 
on the set of objects $f$ of interest.  Suppose that for some (semi-)distance $d_n$ on $\F$, there exist $K>0,\xi>0$, such that for any $\eps>0$ and any $f_1\in\F$ with $d_n(f_0,f_1)\ge \eps$,
\begin{align}
       P_{f_0}^{(n)}\psi_n & \le e^{-Kn\eps^2}   \label{test1}\\
        \sup_{f:\ d_n(f,f_1)\le \xi \eps} P_{f}^{(n)}(1-\psi_n)  \label{test2}
         & \leqa e^{-Kn\eps^2}. 
\end{align}

\begin{theorem}[Thm. 1 in \cite{gvni}] \label{genthm}
Suppose there exists tests $\psi_n$ as in \eqref{test1}-\eqref{test2}, measurable sets $\F_n$ and $C,d>0$, such that, for some $\eps_n\to 0$, $\bar{\eps}_n\to 0$ and $n\eps_n^2\to\pli$, $n\bar{\eps}_n^2\to\pli$,
\begin{align*}
  {\bf (N)} \qquad \quad&  \log N(\bar{\eps}_n,\F_n,d_n)\le d n\bar{\eps}_n^2 \\
  {\bf (S)}  \qquad \quad&  \Pi(\F\backslash\F_n) \le e^{-(C+4)n\eps_n^2}  \\
  {\bf (P)}  \qquad \quad&  \Pi(B_{KL}(f_0,d\eps_n) ) \ge e^{-Cn\eps_n^2} 
\end{align*}
Set $\eps_n^*=\eps_n\vee \bar{\eps}_n$. Then for large enough $M>0$, as $n\to\pli$,
\[
\E_{f_0}\Pi(f: d_n(f,f_0) \leq M\eps_n^* | X^{(n)}) \to 1.
\]
\end{theorem}

\subsection{Applications}

Let us recall that we assume that the compact metric space $\M$ satisfies the conditions of Section \ref{geo} together with the polynomial-type growth \eqref{poly-case} of volume of balls.

\paragraph{Application, Gaussian white noise.}
One observes 
\begin{equation}\label{gwn}
dX^{(n)}(x) = f(x)dx + \frac{1}{\rn} dZ(x),\quad x\in \M.
\end{equation}
 In this case we set $(\mb,\|\cdot\|)=(\L_2,\|\cdot\|_2)$. The prior $W^T$, see \eqref{priorw}, here serves directly as a prior on $f$ (so $w=f$ here).

Here the testing distance is simply $d_n=\|\cdot\|_2$ the $\L_2$-norm in 
$\L_2$. Consider the test
\[
\phi_n=
{\bf 1}_{\{2\int_{\M} \{f_1(x)-f_0(x)\}dX^{(n)}(x) > \|f_1\|^2 - \|f_0\|^2\}}.
\]
Then  \eqref{test1}-\eqref{test2} follow from simple computations.
Also, one can check using Girsanov's formula that for model \eqref{gwn}, the neighborhood  
 $B_{KL}(f_0,\eps_n)$ coincides with an $\L^2$-ball of the same radius. 
%
Recall that the definition of Besov spaces is given in Appendix A.

\begin{theorem}[Gaussian white noise on $(\M,\rho)$] \label{thm-wn}
Let us suppose that $f_0$ is in $B_{2,\infty}^s(\M)$ with $s>0$ and that the prior 
on $f$ is $W^T$ given by \eqref{priorw}. Let $q=1+d/2$ in \eqref{priorT}. Set $\eps_n \sim \bar{\eps}_n\sim (\log{n}/n)^{2s/(2s+d)}$. Then Equations \eqref{priormass}, \eqref{sieve} and \eqref{entropie} are satisfied with  the choice $(\mb,\|\cdot\|_\mb)=(\L^2,\|\cdot\|_2)$. For $M$ large enough, as $n\to+\infty$,
\[ \E_{f_0}\Pi( \|f - f_0\|_2 \ge M\eps_n\ |\ X^{(n)}) \to 0. \]
\end{theorem}

\paragraph{Application, Fixed design regression.}
The observations are
\begin{equation}\label{Fixed design}
Y_i = f(x_i) + \eps_i,\quad 1\le i\le n.
\end{equation}
The design points $\{x_i\}$ are fixed on $\M$ and the variables $\{\eps_i\}$
are assumed to be i.i.d. standard normal. The prior $W^T$, see \eqref{priorw}, here serves directly as a prior on $f$ (so $w=f$ here).

Let us introduce the following 
semi-distance $d_n$. 
For $f_1,f_2$ in $\F$, let us set
\begin{equation}\label{testreg}
d_n(f_1,f_2)^2  =   \int (f_1-f_2)^2 d\P_n^t 
= \frac{1}{n} \sum_{i=1}^{n} (f_1-f_2)^2(x_i).
\end{equation}
Let $\phi_n$ be the likelihood ratio-type test defined by
\[ \phi_n = 1_{\{ \frac{2}{n} \sum_{i=1}^{n} (f_1 - f_0)(x_i) Y_i 
 > \frac{1}{n} \sum_{i=1}^{n} (f_1^2 - f_0^2)(x_i) \} }.\]
This test satisfies, 
 \begin{eqnarray*}
P^{(n)}_{f_0}\phi_n & \leq & G(-\frac{\rn}{2} d_n(\eta_0,\eta_1)  ), \\
\sup_{f\in\F,\ d_n(f,f_1)<d_n(f_0,f_1)/4} P^{(n)}_{f}(1-\phi_n) 
& \leq & G(-\frac{\rn}{4} d_n(\eta_0,\eta_1)  ).\\
\end{eqnarray*}
Also, simple calculations show that for model \eqref{gwn}, the neighborhood 
 $B_{KL}(f_0,\eps_n)$ coincides with an $\L^2(\P_n^t)$-ball of the same radius, which itself contains a sup-norm ball of that radius. 

\begin{theorem}[Fixed design regression on $(\M,\rho)$] \label{thm-regression}
Let us suppose that  $f_0$ is in $B_{\infty,\infty}^s(\M)$, with $s>0$,   and that the prior 
on $f$ is $W^T$ given by \eqref{priorw}. Let $q=1+d/2$ in \eqref{priorT}. Set $\eps_n \sim \bar{\eps}_n\sim (\log{n}/n)^{2s/(2s+d)}$. Then Equations \eqref{priormass}, \eqref{sieve} and \eqref{entropie} are satisfied with $(\mb,\|\cdot\|_\mb)=(\C^0(\M),\|\cdot\|_\infty)$. For $M$ large enough, as $n\to+\infty$,
\[ \E_{f_0}\Pi( d_n(f ,f_0) \ge M\eps_n^*\ |\ X^{(n)}) \to 0. \]
\end{theorem}

\paragraph{Application, Density estimation.}
The observations are a sample 
\begin{equation}\label{density}
(X_i)_{1\le i\le n}\qquad \text{i.i.d.}\quad \sim f,
\end{equation}
for a density $f$ on $\M$. The true density $f_0$ is assumed to be continuous and bounded away from $0$ and infinity on $\M$.  In order to build a prior  on densities, we consider the 
transformation, for any given continuous function $w:\M\to\R$,
\[ f_w^\La(x) := \frac{ \La(w(x)) }{ \int_\M \La(w(u)) d\mu(u) },\quad x\in \M,\] 
where $\La:\R\to (0,\pli)$ is  such that $\log \La$ is Lipschitz on $\R$ and has an inverse $\La^{-1}:(0,\pli)\to\R$. 
For instance, one can take the exponential function as $\Lambda$. Here, the function $w_0$ is taken to be $w_0:=\La^{-1}f_0$.  The prior $W^T$, see \eqref{priorw}, here serves as a prior on $w$'s, which induces a prior on densities via the transformation $f_{w}^\La$. That is, the final prior $\Pi$ on densities 
we consider is $f_{W^T}^{\Lambda}$. In this case we set 
$(\mb,\|\cdot\|)=(\C^0(\M),\|\cdot\|_{\infty})$, the Banach space in which the function $w$ and the prior live.

\begin{itemize}
  \item{ Testing distance.
  It is known from \cite{Birge84} and \cite{LeCam86} that for any two convex sets 
  $\cP_0$ and $\cP_1$ of probability measures, there exist tests $\phi_n$ such that, 
  with $h$ the Hellinger distance,
  \begin{align}\label{testhellinger}
    \sup_{P\in\cP_0} P^n\phi_n & \le \exp(n\log\{1- \frac{1}{2}h^2(\cP_0,\cP_1)\}) \\
    \sup_{P\in\cP_1} P^n(1-\phi_n) & \le \exp(n\log \{1- \frac{1}{2}h^2(\cP_0,\cP_1)\}) 
  \end{align}
  where $h(\cP_0,\cP_1)$ is the infimum of $h(P_0,P_1)$ over $P_0\in\cP_0$ and
  $P_1\in\cP_1$. So, setting $d_n=h$, with the help of the 
  inequality $\log(1-x)\le -x$, we get that \eqref{test1}-\eqref{test2} hold.
    }
  \item{  
 The following is a slight extension of \cite{vvvz09}, Lemma  3.1
\begin{lemma} \label{distco}
For any measurable functions $v,w$, and  a positive function $\La$ on $\R$ such that $\log\La$ is a $L$-Lipschitz function on $\R$, there exists a universal constant $C$ such that
\begin{align*}
h^2(f_v^\La,f_w^\La) & \le L \|v-w\|_{\infty} e^{\|v-w\|_{\infty}/2} \\
K(f_v^\La,f_w^\La) & \le C L \|v-w\|_{\infty} e^{\|v-w\|_{\infty}/2} (1+2L\|v-w\|_{\infty})\\
V(f_v^\La,f_w^\La) & \le C L \|v-w\|_{\infty} e^{\|v-w\|_{\infty}/2} (1+2L\|v-w\|_{\infty})^2.
\end{align*}
\end{lemma}  
With this lemma we see that properties \eqref{entropie}-\eqref{sieve}-\eqref{priormass}
automatically translate up to multiplicative constants into properties {\bf (N)}, {\bf (S)}, {\bf (P)} where $d_n$ is Hellinger's distance.

\begin{proof}
The Hellinger distance between $f_v^\La$ and  $f_w^\La$ can be written
\[h(f_v^\La,f_w^\La) = \| \frac{\sqrt{ \La(v) }}{\| {\sqrt{ \La(v) }}\|_2} - 
\frac{\sqrt{ \La(w) }}{\| {\sqrt{ \La(w) }}\|_2} \|_2
\le 2  \frac{\| \sqrt{ \La(v) } - \sqrt{ \La(w) }\|_2} {\| {\sqrt{ \La(w) }}\|_2}. \]
The inequality $e^x\le 1+x e^x$ for $x\ge 0$ implies
\begin{align*}
 \sqrt{ \La(v) } - \sqrt{ \La(w) } & = \sqrt{ \La(w) } \left[e^{(\log\La(v) - 
 \log\La(w) )/2} - 1\right]\\
 & \le \sqrt{ \La(w) } \left[ e^{L\|v-w\|_{\infty}/2} - 1\right]\\
 & \le \frac{L}{2} \sqrt{ \La(w) }  \|v-w\|_{\infty} e^{\|v-w\|_{\infty}/2}. 
\end{align*}
This leads to the first inequality. By Lemma 8 in \cite{gvv07}, for any densities $p,q$
\[ K(p,q) \leqa h(p,q) \left(1+ \log\|\frac{p}{q}\|_{\infty}\right),\qquad V(p,q) \leqa h(p,q) \left(1+  \log\|\frac{p}{q}\|_{\infty}\right)^2. \]
Since $\log \La$ is $L$-Lipschitz, we have that
\[e^{-L\|v-w\|_{\infty}} \int \La(w) d\mu  \le \int \La(v) d\mu \le
e^{L\|v-w\|_{\infty}} \int \La(w) d\mu.\]
Inserting this into the following chain of inequalities,
\begin{align*}
\log \|\frac{f_v^\La}{f_w^\La} \|_{\infty} & \le \|\log \frac{f_v^\La}{f_w^\La} \|_{\infty} \\
& \le \| \La(v)-\La(w) - \log\frac{\int \La(v)d\mu}{\int \La(w)d\mu} \|_{\infty}\\
& \le 2L\|v-w\|_{\infty}.
\end{align*}
This proves the Lemma. \qed
\end{proof}

   }
 \end{itemize}

\begin{theorem}[Density estimation on $(\M,\rho)$] \label{thm-density}
Suppose the true density $f_0$ is a continuous function bounded away from $0$ and $\pli$ over $\M$. Let the prior $\Pi$ be the law on densities induced by $f_{W^T}^\Lambda$, with $W^T$ defined by \eqref{priorw}. Let $q=1+d/2$ in \eqref{priorT} and $\La$ be a positive invertible function  with $\log\La$ Lipschitz on $\R$. Let us suppose that   $\La^{-1}f_0$ is in $B_{\infty,\infty}^s(\M)$, $s>0$. Set $\eps_n\sim\bar{\eps}_n
\sim (\log{n}/n)^{s/(2s+d)}$. Then Equations \eqref{priormass}, \eqref{sieve} and \eqref{entropie} are satisfied with 
$(\mb,\|\cdot\|_\mb)=(\C^0(\M),\|\cdot\|_\infty)$. For $M$ large enough, as $n\to+\infty$,
\[ \E_{f_0}\Pi( h(f ,f_0) \ge M\eps_n^*\ |\ X^{(n)}) \to 0. \]
\end{theorem}

\paragraph{Proofs.} The proofs of Theorems \ref{thm-wn} and \ref{thm-regression} directly follow from the results in Section \ref{section_gc}. Indeed, for the white noise case, apply Theorem \ref{genthm} with $\F=\L^2$, $d_n=\|\cdot\|_2$, $\F_n$ the set $B_n$ defined in Section \ref{section_gc} and $f_0=w_0\in B_{2,\infty}^{s}$. The regression case is similar with $\F=\C^0(\M)$ equipped with the sup-norm and $f_0=w_0\in B_{\infty,\infty}^{s}$.

In the density case, the key property is the following inclusion, 
for some $c>0$ and any small enough $\eps>0$ (recall that $w_0=\La^{-1}f_0$ so that
$f_0=f_{w_0}^{\La}$)
\begin{align}\label{distcomp}
\{f_{w}^{\La}:\ \|w - w_0\|_{\mb} \le \eps \}  \quad \subset & \quad\{f:\ d_n(f,f_0)\le c\eps\} \cap
 B_{KL}(f_0,c\eps).
\end{align}
This means that the testing distance $d_n$ on densities and the KL-divergence properly relate to the Banach space distance $\|\cdot\|_{\mb}$ on the set of $w$'s.
Now the inclusion property \eqref{distcomp} is clear in view of the definition of $B_{KL}$ and Lemma \ref{distco}. Therefore, the inequalities \eqref{priormass}-\eqref{sieve}-\eqref{entropie} of Section \ref{section_gc} related to functions $w$  automatically translate into the properties {\bf (N)}, {\bf(S)}, {\bf (P)} for functions $f$ needed for Theorem \ref{genthm}. The remainder of the proof is as for the white noise and regression case.
\qed 

In the case that $\M$ is a compact connected orientable manifold without boundary, minimax rates of convergence have been obtained in \cite{efromovich}, where Sobolev balls of smoothness index $s$ are considered and data are generated from a regression setting. In particular, in this framework, our procedure is adaptive in the minimax sense for Besov regularities, up to a logarithmic factor.
 
We have obtained convergence rates for the posterior distribution associated to the geometrical prior in a variety of statistical frameworks. Obtaining these rates does not presuppose any a priori knowledge of the regularity of the function $f_0$. Therefore our procedure is not only nearly minimax, but also nearly adaptive.

  Note also that another attractive property of the method is that it does not assume a priori any (upper or lower) bound on the regularity index $s>0$. This is related to the fact that approximation is via the spaces $\bH_t$, which are made of (super)-smooth  functions.

\section{Lower bound for the rate} \label{sec-lb}

Works obtaining (nearly-)adaptive rates of convergence for posterior distributions are relatively recent and so far were obtained for density or regression on subsets of the real line or the Euclidian space. Often, logarithmic factors are reported in the (upper-bound) rates, but it is unclear whether the rate must include such a logarithmic term.  We aim at giving an answer to this question in our setting by providing a lower bound for the rate of convergence of our general procedure. This lower bound implies that the rates obtained in Section \ref{section_rates} are, in fact, sharp. One can conjecture that the same phenomenon appears for hierarchical Bayesian procedures with randomly rescaled Gaussian priors when the initial Gaussian prior has a RKHS which is made of super-smooth functions (e.g. infinitely differentiable functions), for instance the priors considered in \cite{vvvz09}, \cite{VincentJudith-np}.

For simplicity we consider the Gaussian white noise model
\begin{equation*}
dX^{(n)}(x) = f(x)dx + \frac{1}{\rn} dZ(x),\quad x\in \M.
\end{equation*}
We set $(\mb,\|\cdot\|)=(\L_2(\M),\|\cdot\|_2)$. As before, for this model the prior sits on the same space as the function $f$ to be estimated, so $w=f$. 


\begin{theorem}[Gaussian white noise on $(\M,\rho)$, lower bound] \label{thm-lb}
Let   $\eps_n =(\log{n}/n)^{s/(2d+s)}$ for $s>0$ and let the prior on $f$ be the law induced by $W^T$, see \eqref{priorw},  with $q>0$ in \eqref{priorT}. Then there exist $f_0$  in the unit ball of  $B_{2,\infty}^s(\M)$ and a constant $c>0$
such  that  
\[ \E_{f_0}\Pi( \|f - f_0\|_2 \le c(\log{n})^{0\vee (q-1-\frac{d}{2})}\eps_n\ |\ X^{(n)}) \to 0. \]
\end{theorem}

As a consequence, for any prior of the type \eqref{priorw} with any $q>0$ in \eqref{priorT}, the posterior convergence rate cannot be faster than $\eps_n$ above. If $q$ is larger than 
$1+d/2$, the rate becomes even slower than $\eps_n$.

\begin{rem}
More generally, an adaptation of the proof of Theorem \ref{thm-lb} yields that, for any `reasonable' prior on $T$, in that, for $\eps_n \sim (\log{n}/n)^{s/(2d+s)}$, it holds  
\[ \Pi[\|f-f_0\|_2\le \eps_n ] \ge e^{-Cn\eps_n^2},\] 
then $\Pi[\|f-f_0\|_2\le c\eps_n \given X]\to 0$ for small enough $c>0$.  This condition is the standard `prior mass' condition in checking upper-bound rates, see \eqref{priormass}. Note that the previous display is automatically implied if the prior satisfies $\Pi[\|f-f_0\|_2\le \eps_n^* ] \ge e^{-Cn{\eps_n^*}^2}$ for $\eps_n^*=n^{-s/(2d+s)}$, or more generally for any rate at least as fast as $\eps_n$. For instance, this can be used to check that taking a uniform prior on $(0,1)$ as law for $T$ leads to the same lower bound rate. 
\end{rem}

\begin{proof} 
We use a general approach to prove lower bounds for posterior measures introduced in \cite{ic08} (see \cite{ic08}, \cite{spa} for examples). The idea is to apply the following lemma (Lemma 1 in \cite{gvni}) to the sets 
$\{f\in\mb,\ \|f-f_0\|_\mb\le \zeta_n\}$, for some rate $\zeta_n\to 0$ and $f_0$ in $B_{2,\infty}^s$, with $s>0$.
\begin{lemma} \label{lemba}
If $\al_n\to 0$ and $n\al_n^2\to \pli$  and if $B_n$ is a measurable set such that
\[\Pi(B_n)/\Pi(B_{KL}(f_0,\al_n)) \leq e^{-2n\al_n^2},\]
then $\E_{f_0}\Pi(B_n\ | \ X^{(n)}) \to 0$ as $n\to\pli$.
\end{lemma}
In our context this specializes as follows. Let  $\al_n\to 0$ and $n\al_n^2\to \pli$. Suppose that, as $n\to\pli$, 
\[ \frac{\Pi(\| f-f_0 \|_2\le \zeta_n) }{ \Pi(\|f-f_0\|_2\le \al_n)} = o(e^{-2 n\al_n^2}). \]
Then $\zeta_n$ is a lower bound for the rate of the posterior in that, as $n\to\pli$, 
\[ \E_{f_0}\Pi(\| f-f_0 \|_2\le \zeta_n) \ | \ X^{(n)}) \to 0.\]

We first deal with the case where $q\le 1+d/2$. In this case let us choose $\al_n=2\eps_n$, where $\eps_n=(\log n/n)^{2s/(2s+d)}$. In Section \ref{section_gc}, we have established in \eqref{pmcond} that, for the prior $W^T$ with $q\le 1+d/2$ in  \eqref{priorT}, there exists $C>0$ with 
\[ \Pi(\|f-f_0\|_2 \le \eps_n) = \mathbb{P}(\|W^T-w_0\|_{\mb}\le \eps_n)
 \ge e^{-Cn\eps_n^2}.\]
So it is enough to 
show that, for some well-chosen $\zeta_n\to 0$, 
\begin{equation} \label{lb-goal} 
\Pi(\| f-f_0 \|_2\le \zeta_n)) =o(e^{-(8+C)n\eps_n^2}). 
\end{equation}
We would like to take $\zeta_n=c\eps_n$, for some (small) constant $c>0$.
In order to bound from above the previous probability, we write  
\begin{align*}
\Pi[ \| f-f_0 \|_2\le \zeta_n ] & = \int_0^1  \Pi[ \| W^t-f_0 \|_2\le \zeta_n ] g(t)dt  \\
& \le \int_0^1  \exp\left[ -\phi^t_{f_0}(\zeta_n)  \right] g(t)dt.
\end{align*}
We separate the above integral in two parts. The first one is $\mathcal{T}_1:=\{\mu_n\le t\le B t^*_n\}$, where $t_n^*$ is a similar cut-off as in the upper-bound proof $t_n^*=\zeta_n^{2/s}\log(1/\zeta_n)$. On $\mathcal{T}_1$, one can bound from below $\phi^t_{f_0}(\zeta_n)$ by its small ball probability part $\phi_0^t(\zeta_n)$. 
Moreover, thanks to relation (3.1) in \cite{KuelbsLi93}, we have, for any $\la>0$ and $t\in(0,1]$, 
\[ \phi_0^t(\zeta_n) = - \log\P[\|W^t\|_2<\zeta_n] \ge 
H(2\zeta_n,\la\bH_t^1,\|\cdot\|_2)-\frac{\la^2}{2}.\]  
Set $\la=1$ and recall from Remark \ref{rementropy} that the lower bound on the entropy can be used for any $t$ regardless of the value of $\eps$. This yields, for large enough $n$, if $\zeta_n=o(1)$,
 \[ \phi_0^t(\zeta_n) \geq C (B t_n^*)^{-d/2} \log^{1+d/2}(1/\zeta_n) - \frac{1}{2}
 \geq C B^{-d/2}\zeta_n^{-d/s}\log(1/\zeta_n). \]
Thus we obtain
\begin{align*}
\int_{0}^{B t_n^*} \exp\left[ -\phi^t_{f_0}(\zeta_n)  \right] g(t)dt
& \le e^{-C B^{-d/2}\zeta_n^{-d/s}\log(1/\zeta_n) } \int_{0}^{B t_n^*} g(t)dt\\
& \le  e^{- C B^{-d/2}\zeta_n^{-d/s}\log(1/\zeta_n) }.  
\end{align*}  
This is less than $e^{-(8+C)n\eps_n^2}$ provided $\zeta_n=\kappa\eps_n$ and $\kappa>0$ is small enough.

It remains to bound the integral from above on 
$\cT_2:=\{Bt_n^* \le t \le 1\}$.  Here we bound $ \phi^t_{f_0}(\zeta_n)$ from below by its approximation part. For any $t\in \cT_2$, 
\[ \phi^t_{f_0}(\zeta_n) \ge \frac{1}{2} \cdot
\inf_{h\in \bH_t,\ \|h-f_0\|_2<\zeta_n} \|h\|_{\bH_t}^2. \]
 We prove in Appendix B, see Theorem \ref{LowerboundA}, that there exist constants $c,\; C$ and $f_0$ in $B_{2,\infty}^s(\M)$ such that 
 \begin{equation} \label{lb-approx}
   \phi^t_{f_0}(\zeta_n) \ge C\zeta_n^2 e^{c\zeta_n^{-2/s} t}.
 \end{equation}
 
Now, under \eqref{lb-approx} for the previous fixed function $f_0$, taking $\zeta_n=\kappa \eps_n$ for small (but fixed) enough  $\kappa$, it holds, when $t$ belongs to $\cT_2$, 
 \[\phi^t_{f_0}(\zeta_n)\ge C(\kappa \eps_n)^a 
e^{c \kappa^{-2/s} B \log(1/\zeta_n)}.\] For $\kappa$ small enough, this is larger that any given power of $\eps_n$. In particular, it is larger than $(8+C)n\eps_n^2$ if the (upper-bound) rate $\eps_n$ is no more than polynomial in $n$, which is the case here since $\eps_n=(\log n/n)^{s/(2s+p)}$.
We have verified that \eqref{lb-goal} is satisfied, which gives the desired lower bound result when $q\le 1+d/2$  using Lemma \ref{lemba}. 

In the case that $q> 1+d/2$, the proof is the similar, except that the exponent of the logarithmic factor in \eqref{pmcond} has now  the  power $q-d/2$, due to the assumption on the prior density $g$, and that  $\eps_n$ is now replaced by 
$\tilde{\eps}_n=(\log{n})^{q-1-\frac{d}{2}}\eps_n$.
\end{proof}

\section{Appendix A: Besov spaces} \label{appendix-besov}

We follow the paper \cite{CKP} to introduce the Besov spaces $B^s_{pq}$ in this setting
with $s>0$, $1\le p \le \infty$ and $0<q\le \infty$. 
To do so, let us introduce a (Littlewood-Paley) function
$\Phi\in C^\infty(\mathbb{R}_+)$ such that
\begin{align}
&\supp \Phi \subset   [0, 2] ,\;\;
\Phi\ge 0  \hbox{ for } \nu\ge 1,\;\;
|\Phi(\lam)| =1 \;\hbox{ for } \lam\in [0, 1], \label{cond_phi}
\end{align}
Set $\Phi_j(\lambda):= \Phi(2^{-j}\lambda)$ for $j\ge 1$.


\begin{definition}\label{def-B-spaces}
Let $s>0$, $1\le p\le \infty$, and $0<q \le \infty$.
The Besov space  $B_{pq}^{s}=B_{pq}^{s}(L)$
is defined as the set of all $f \in \L^p(\M,\mu)$ such that
\begin{equation}\label{def-Besov-space1}
\|f\|_{B_{pq}^{s}} :=
\Big(\sum_{j\ge 0} \Big(2^{s j}
\| \Phi_j(\sqrt{L}) f(\cdot) - f(\cdot) \|_{\L_p}
\Big)^q\Big)^{1/q} <\infty.
\end{equation}
Here the $\ell^q$-norm is replaced by the sup-norm if $q=\infty$.
\end{definition}

\begin{rem}
This definition is independent of the choice of $\Phi$. Actually if
%
$\cE_t(f)_p$ denotes the best approximation of
$f \in \L^p$ from $\Sigma_t$, that is,
$$\cE_t(f)_p:=\inf_{g\in\Sigma_t}\|f-g\|_p.$$
(here $\L^\infty$ is identified as the space $\UCB$ of all uniformly continuous
and bounded functions on $M$)
%
%
%
then
$$B_{pq}^s:=\{ f\in \; \L_p\; / \|f\|_{A_{pq}^s}:= \|f\|_p +
\Big(\sum_{j\ge 0}\big(2^{s j}\cE_{2^j}(f)_p \big)^q\Big)^{1/q}<\infty \}.$$
\end{rem}


\section{Appendix B: Entropy properties} \label{appendix-entropy}
\subsection{Covering number, entropy, $\epsilon-$net.}

Let $(X, \rho)$ be a metric space. For $\epsilon >0$ the covering number $N(\epsilon, X) $ is the smallest number of balls of radius
$\epsilon $ covering $X$. The entropy $H(\epsilon, X)$ is by definition $H(\epsilon, X)= \log_2 N(\epsilon, X) .$ 
\\
An important result of this section paper is the link between the covering number $N(\epsilon, \M,\rho) $ of the space $\M$,   and 
$H(\epsilon, B, \L^p )$ the entropy number of the unit ball $B$ of some functional space, computed in the $\L^p$ metric.
\\
An $\epsilon-$net $\Lambda \subset X$ is a set such that $\xi \neq \xi', \; \xi, \xi' \in \Lambda $ implies $\rho(\xi, \xi') >\epsilon.$
A maximal  $\epsilon-$net $\Lambda$, is a an $\epsilon-$net such that $\hbox{ for all } x \in X\setminus \Lambda, \; \Lambda \cup \{x\}$
is no more an $\epsilon-$net.
So, for a maximal $\epsilon-$ net, $\Lambda$ we have :
$$ X \subset \cup_{\xi \in \Lambda} B(\xi, \epsilon) , \; \xi \neq \xi', \; \xi, \xi' \in \Lambda \Rightarrow B(\xi, \epsilon/2) \cap B(\xi', \epsilon/2)
=\emptyset.$$
Hence, for $\Lambda_\epsilon$    a maximal $\epsilon-$ net  we have : 
$$N(\epsilon/2, X) \geq card(\Lambda_\epsilon) \geq N(\epsilon, X).$$ 
Now if $(X, \rho)$ is a doubling metric space then we have the following property :
If $x_1, \ldots,x_N \in B(x,r)$ are such that, $ \rho(x_i,x_j) > r 2^{-l}$  ($l \in \N$)  clearly $B(x,r) \subset  B(x_i, 2r)=B(x_i,  2^{l+2} (r 2^{-l-1}))$
and the balls $B(x_i ,  r 2^{-l-1})$ are disjoint and contained in $B(x, 2r)$. so :
\begin{equation}\label{covc}
 N 2^{-(l+2)D}  |B(x,r)| \leq \sum_{i=1}^N |B(x_i, r2^{-l-1})| \leq B(x, 2r)| \leq 2^D |B(x,r)|
 \end{equation}
If $\Lambda_{r2^{-l}}$ is any $r2^{-l}-$net then :
$Card(\Lambda_{r2^{-l}}) \leq  2^{(l+3)d}N(X,r).$
So if $\Lambda_{\epsilon}$ is any maximal $ \epsilon-$net  and for $ l\in \N, \;
 \Lambda_{2^l\epsilon}$ is any maximal $2^l \epsilon-$net   then :
\begin{equation}\label{entr}
N(X, \epsilon 2^l) \leq 
N(X, \epsilon) \leq  Card (\Lambda_{ \epsilon}) \leq  2^{(l+3)D}N(X, 2^l\epsilon) \leq 2^{(l+3)D}Card (\Lambda_{2^l\epsilon}).
 \end{equation}
For $l=0$
$$ 2^{-3D}Card (\Lambda_{ \epsilon}) \leq  N(X, \epsilon) \leq Card (\Lambda_{ \epsilon}) .$$
So for any $ \epsilon >0, $ and for any  maximal $ \epsilon-$net $\Lambda_{ \epsilon}, \;  Card (\Lambda_{ \epsilon}) $ and $  N(X, \epsilon) $ are of the same order.
\\
Moreover clearly, taking $r=1$ in (\ref{covc}), so that $B(x,1)=\M,$ we get:
\begin{equation}\label{nbr}
N(\delta,\M) \leq 4^D(\frac 1{\delta})^D
\end{equation}
\begin{rem}\label{reg}
In number of examples (for instance for compact Riemannian manifolds) there exist absolute constants :
$  c_1>, c_2>0,  d>0 $ such that 
$$\hbox{ for all } x \in \M, \;  \hbox{ for all } \;  0 <r  \leq 1, \; c_1 r^d \leq | B(x,r)| \leq c_2 r^d.$$
Necessarily, $ d \leq D$ since using \eref{inv} $|B(x,r)| \geq (r/2)^D$; so $c_1 r^d \geq (r/2)^D, \; 0<r \leq 1,$.
\\
Let  $(B(x_i,\epsilon))_{i\in I}$ be a minimal covering of $\M$; we have
$$ 1= |\M| \leq \sum_{i\in I}  |B(x_i, \epsilon)|  \leq N(\epsilon, \M) c_2 \epsilon^d.$$
Now if  $\Lambda_{\epsilon}$ is any maximal $ \epsilon-$net, we have :
$$ 1= |\M| \geq \sum_{\xi \in \Lambda_{\epsilon}} B(\xi, \epsilon/2)  \geq card(\Lambda_{\epsilon}) c_1 (\epsilon/2)^d$$
As a conclusion, in the case \eref{poly-case}:
\begin{equation}\label{poly-entropie} \frac{1}{c_2}  (\frac 1\epsilon)^d  \leq N(\epsilon, \M) \leq  card(\Lambda_{\epsilon})  \leq  \frac{2^n}{c_1} (\frac 1\epsilon)^d.
\end{equation}
\end{rem}

\subsection{Dimension of spectral spaces, covering number, and trace of $P_t$.}
Let us now use the  'heat kernel' assumptions. The  following proposition gives the link between  the covering number $ N(\delta, \M)$ of the underlying space $\M$, the behavior of
the trace of $e^{-tL}$ and the dimension of the spectral spaces.
Let us define:
$$ \Sigma_\lambda = \otimes_{\sqrt{\lambda_k} \leq \lambda} \cH_{\lambda_k}.$$
Clearly the projector  $P_{\Sigma_\lambda}$ is a kernel operator and 
$$ P_{\Sigma_\lambda}(x,y) = \sum_{\sqrt{\lambda_k} \leq \lambda} P_k(x,y)$$
Then one can prove the following bounds (see \cite{CKP}, Lemma 3.19): For  any $\lambda \geq 1, $ and $\delta= \frac 1\lambda,$
\begin{equation}\label{proj}
 \exists C_2, C'_2, \quad \hbox{such that}\quad
\frac{C'_2}{|B(x,\delta)|}  \leq P_{\Sigma_\lambda}(x,x) \leq \frac{C_2}{|B(x,\delta)|}
\end{equation}

Let us recall that $ Tr(e^{-tL}) =\sum_k e^{-\lambda_kt} dim(\cH_k) .$ In addition we have $\int_{\M} P_{t}(x, x) d\mu(x) = Tr(e^{-tL})$.
Moreover , as 
$$P_t(x,x) = \int_{\M} P_{t/2}(x,u) P_{t/2}(u,x) d\mu(u)= \int_{\M}  (P_{t/2}(x,u) )^2 d\mu(u)$$ we have : 
$$Tr(e^{-tL})
= \int_{\M} P_t(x,x)  d\mu(x) = \int_{\M} \int_{\M}   (P_{t/2}(x,u) )^2 d\mu(u)d\mu(x) =  \|e^{- \frac t2 L}\|_{HS}^2. $$
where $\| \; \|_{HS}$ stands for the Hilbert-Schmidt norm.
\begin{proposition}\label{prop5}
\begin{enumerate}
\item For $ \lam \geq 1, \quad \delta = \frac 1\lambda, \quad$ 
\begin{equation}\label{finite-dim-2}
  C'_2  \int_{\M}  \frac 1{|B(x, \delta)|} d\mu(x) \leq \dim (\Sigma_\lam) =   \int_{\M}  P_{\Sigma_\lambda}(x,x)  d\mu(x)  \leq C_2  \int_{\M}  \frac 1{|B(x, \delta)|} d\mu(x)
\end{equation}
\item
\begin{equation}\label{fin-dim-1}
2^{-2D} N(\delta, \M) \leq 
2^{-2D} card(\Lambda_\delta) \leq  \int_{\M}  \frac 1{|B(x, \delta)|} d\mu(x) \leq 2^D card(\Lambda_\delta)
\leq 2^{4D} N(\delta, \M) 
\end{equation}
where $\Lambda_\delta$ is any $\delta-$maximal net.
\item
$$ C'_1 \int_{\M} \frac 1{|B(x,\sqrt t)|} d\mu(x) \leq Tr(e^{-tL}) \leq C_1
 \int_{\M} \frac 1{|B(x,\sqrt t)|} d\mu(x)
$$

\end{enumerate}
\end{proposition}

{\bf Proof of the Proposition:}
1.  is a consequence of (\ref{proj}).
Let us prove 2. : \\
Let  $\Lambda_\delta$ be any $\delta-$maximal net.
$$\sum_{\xi \in \Lambda_\delta }  \int_{B(\xi, \delta/2)}  \frac 1{|B(x, \delta)|} d\mu(x) \leq 
 \int_M  \frac 1{|B(x, \delta)|} d\mu(x) \leq \sum_{\xi \in \Lambda_\delta }  \int_{B(\xi, \delta)}  \frac 1{|B(x, \delta)|} d\mu(x) $$
But :
$$ x \in B(\xi, \delta/2) \Longrightarrow  B(x, \delta) \subset B(\xi, 2\delta), \quad 
\hbox{so} \quad 
 \frac 1{|B(x, \delta)|} \geq  \frac{2^{-2D}}{|B(\xi, \delta/2)|}$$
and in the same way :
$$ x \in B(\xi, \delta) \Longrightarrow  B(\xi , \delta) \subset B(x, 2\delta), \quad
\hbox{so} \quad 
 \frac 1{|B(x, \delta)|} \leq  \frac{2^{D}}{|B(\xi, \delta)|}$$
This implies :
$$2^{-2D} card(\Lambda_\delta)  \leq \int_{\M}  \frac 1{|B(x, \delta)|} d\mu(x) \leq 2^D card(\Lambda_\delta)$$
3.  is a consequence of (\ref{hker}).

\qed

The former results can be summarized in the following corollary:

\begin{corollary}\label{equi}
$$ Trace (e^{- \delta^2 L} ) \sim dim( \Sigma_\lambda) \sim N(\delta,\M) ; \quad \delta = \frac 1\lambda$$

\end{corollary}

%
%
%
%
%
%
%
%
%
%
%
%
%
%
\subsection{Connection between covering number of $\M$ and entropy of $\bH^1_t$.}
In  this section we establish the link between the covering number $N(\epsilon, \M) $ of the space $\M$,   and 
$H(\epsilon, \bH^1_t, \L^p )$ for $p=2,\; \infty$ stated in Theorem \ref{entropyconnection}, which we recall here.
\\
 Let us suppose that for some $\mu >0, \; a>0, \;   \epsilon^\mu \leq at.$
Then there exists  $\epsilon_0 >0,$ such that for all $ 0<\epsilon\leq \epsilon_0, \quad$
$$ H(\epsilon, \bH^1_t, \L^2) \sim   H(\epsilon, \bH^1_t, \L^\infty) \sim    N(\delta(t,\epsilon), \M)  \cdot \log \frac 1\epsilon
\quad \hbox{where}\quad \frac 1{\delta(t,\epsilon)}=
  \sqrt{\frac 1t  \log (\frac 1\epsilon)}.
$$

\noindent Notice, of course, that one can replace $N(\delta(t,\epsilon), \M) $ at any place by $ card(\Lambda_{\delta(t,\epsilon) }) ,$ where $ \Lambda_{\delta(t,\epsilon) } $ is a   maximal $\delta(t,\epsilon)-$net.
Also, since $\mu(\M)=1$, we have  
$$H(\epsilon, \bH^1_t, \L^2)   \leq H(\epsilon, \bH^1_t, \L^\infty).$$
 So the proof will be done in two steps:\\
 We prove the lower bound for $H(\epsilon, \bH^1_t, \L^2) $ in the next subsection, using
 Carl's inequality.
\\
We prove next the upper bound for $H(\epsilon, \bH^1_t, \L^\infty).$

\subsubsection{ Proof of the theorem: Lower estimates for  $ H(\epsilon, \bH^1_t, \L^2)  $. }

Let us recall some classical facts: see the following references \cite{Carl1981}, \cite{Carlstephani}.
For any  subset $X$  of a metric space, we define , for any $k \in \N$
$$ e_k(X) = \inf \{\epsilon \geq 0, \; \exists \;  2^k  \quad \hbox{balls of radius } \; \epsilon, \quad \hbox{covering}\  X.\} .$$
Clearly 
$$\epsilon <e_k(X) \Longrightarrow H(\epsilon, X) >k$$
Now for the special case of a compact positive selfadjoint  operator $T : \bH  \mapsto  \bH$   we have the following  Carl (cf \cite{Carl1981} ) inequality relating
$e_k(T(B))$ where $B$ is the unit ball of $\bH$ and the eigenvalues  $ 0 \leq \mu_1 \leq \mu_2,\ldots$ (possibly repeated with their multiplicity order) of $T$:
\begin{equation}\label{carl}
 \hbox{ for all } k \in \N^*, \; n \in \N^*, \quad e_k(T(B)) \geq 2^{-\frac k{2n}} \prod_{i=1}^n \mu_i^{1/n}
 \end{equation}
 In our case, let us take : $T= P_{t/2}, \; \mu_i =e^{-t/2 \lambda_i } , T(B)= \bH_t^1$. Let us fix :
 
$$\lambda = \sqrt{  \frac 1t \log \frac 1\epsilon} = \frac 1\delta = \frac 1{                                                                                                             \delta(t,\epsilon)}$$
$$ n= dim ( \Sigma_\lambda ); \quad
 k  \sim n \log \frac 1\epsilon \frac 1{\log 2}     $$
%
 Carl's inequality gives:
 $$ e_k \geq  2^{-\frac k{2n}} e^{- \frac 1n \sum_{ t\lambda_i \leq \log \frac 1{\sqrt \epsilon}} t/2 \lambda_i}  \geq \epsilon $$
So
$$ H(\epsilon, \bH_t^1, \L^2) \geq k  \sim n \log \frac 1\epsilon \frac 1{\log 2}  \sim dim ( \Sigma_\lambda )  \log \frac 1\epsilon   $$
 but  by  (\ref{equi} ) $ dim ( \Sigma_\lambda ) \sim N( \delta, \M), \quad \delta = \frac 1\lambda .$ 
 So :
$$ H(\epsilon, \bH_t^1, \L^2) \gtrsim  \log \frac 1\epsilon N( \delta, \M),
 \quad  \frac 1\delta =
\lambda = \sqrt{  \frac 1t \log \frac 1\epsilon}. $$

%
 
\subsubsection{Proof of the theorem: Upper estimate for  $H(\epsilon,H^1_t,\L_\infty).$}
\label{proof-entropy} 
 We recall the notations introduced in  Section \ref{geo} (especially \ref{smoothoperator}).
  Let us suppose : $ \epsilon^\mu  \leq at, \mu > 0  , a>0.$\\
   First, we prove that  for all $ \epsilon >0,$ small enough, there exists $  \delta\; ( \sim \delta(t, \epsilon) := \sqrt{\frac 1t \log \frac 1\epsilon}) $  such that  
  $$\hbox{ for all } f \in \bH^1_t, \;  \| \Phi(\delta \sqrt L)f- f \|_\infty \leq \frac \epsilon2.$$
  
  In a second step we use \eref{wav} to expand on the $|B(\xi, \delta)| D^\delta_\xi$'s:
  $$ \Phi(\delta \sqrt L)f(x) = \sum_{\xi \in  \Lambda_{\gamma \delta}} \Phi(\delta \sqrt L)f(\xi) |B(\xi, \delta)| D^\delta_\xi(x) .$$
  
  In a third step, we use a family of points of $\Sigma_{\frac 1\delta}$ as centers of balls of radius $\epsilon/2$ covering 
  $ \Phi(\delta \sqrt L) ( \bH^1_t)$ so that the balls centered in these points is an $\epsilon-$ covering in $\L^\infty$ norm of
 $ \bH^1_t.$

  The next lemma gives  evaluations of $ \| \Phi(\delta \sqrt L)f- f \|_\infty $ and $\|\Phi(\delta \sqrt L) (  \bH^1_t )\|_\infty.$ 
  
\begin{lemma}\label{lemm}

$ \hbox{ for all }  f \in \bH^1_t$
\begin{enumerate}
\item
$$  \|\Phi(\delta \sqrt L) f  \|_\infty \lesssim   \frac 1{t^{D/4} } $$
\item
$$  \|\Psi(\delta \sqrt L) f   \|_\infty \lesssim    e^{-  \frac t{8\delta^2}}\frac 1{\delta^{D/2}}$$
\item
$$  \| \Phi(\delta \sqrt L)f- f \|_\infty  \leq \sum_{j\geq 0} \|\Psi( 2^{-j}\delta \sqrt L) f \|_\infty \lesssim   \frac 1{\delta^{D/2}} e^{-\frac A4} A^{-1}, \quad A=  \frac{ t}{8\delta^2}$$

\end{enumerate}

\end{lemma}
\noindent
{\bf Proof of the lemma:}
\\
First,
$ f \in \bH^1_t$ so $ \; f=\sum_k \sum_l a_k^l e^l_k(.)e^{-\lambda_kt/2}, \;  \sum_k \sum_l  | a_k^l |^2 \leq 1$.
As $\Phi(\delta \sqrt L)(x,y) =\sum_k \Phi(\delta \sqrt{\lambda_k}) P_k(x,y)$,
$$\Phi(\delta \sqrt L) f(x)= \langle \Phi(\delta \sqrt L)(x,.), f(.) \rangle =\sum_k \sum_l \Phi(\delta \sqrt{\lambda_k})  a_k^l e^l_k(x)e^{-\lambda_kt/2},\hbox{ hence}$$
\begin{align*}  |\Phi(\delta \sqrt L) f(x)| &\leq  (\sum_k \sum_l |a_k^l |^2)^{1/2}
(\sum_k  e^{-\lambda_kt} \Phi^2(\delta \sqrt{\lambda_k})  \sum_l (e^l_k(x)^2 )^{1/2}\\
&\leq (\sum_k  e^{-\lambda_kt} \Phi^2(\delta \sqrt{\lambda_k}) P_k(x,x) )^{1/2}\\
&\leq   (\sum_k   \Phi^2(\delta \sqrt{\lambda_k}) P_k(x,x) )^{1/2} \wedge  (\sum_k  e^{-\lambda_kt}  P_k(x,x) )^{1/2}
\\
&=  [ \Phi^2(\delta \sqrt L) (x,x)  \wedge  P_t(x,x) ]^{1/2}\\
&\leq  \frac{ \sqrt{C(\Phi^2)}}{\sqrt{|B(x, \delta)|}} \wedge \frac{ \sqrt{C_1}}{\sqrt{|B(x, \sqrt t)|}} \lesssim  \frac{1}{t^{D/4}}
\end{align*}
 using \eref{hker}, \eref{FC} and (\ref{inv}).

In the same way :
$$\Psi(\delta \sqrt L) f(x)= \langle \Psi(\delta \sqrt L)(x,.), f(.) \rangle =\sum_k \sum_l \Psi(\delta \sqrt{\lambda_k})  a_k^l e^l_k(x)e^{-\lambda_kt/2},\hbox{ hence}$$
\begin{align*} |\Psi(\delta \sqrt L) f(x)| &\leq  (\sum_k \sum_l |a_k^l |^2)^{1/2}
(\sum_k  e^{-\lambda_kt  } \;  \Psi^2(\delta \sqrt{\lambda_k})  \sum_l (e^l_k(\xi))^2 )^{1/2}\\
&\leq   e^{- \frac 1{4\delta^2} t/2}(\sum_k   \Psi^2(\delta \sqrt{\lambda_k}) P_k(x,x) )^{1/2}
\\
&=e^{-  \frac t{8\delta^2}}  [ \Psi^2(\delta \sqrt L) (x,x)]^{1/2}\\
&\leq C(\Psi^2)  e^{-  \frac t{8\delta^2}} \frac 1{|B(x,\delta)|^{1/2}} 
\\ 
&\leq 2^{D/2}C(\Psi^2)   e^{-  \frac t{8\delta^2}}\frac 1{\delta^{D/2}}
\\ 
& \lesssim e^{-  \frac t{8\delta^2}}\frac 1{\delta^{D/2}}.
 \end{align*}

So
$$\sum_{j\geq 0} \|\Psi( 2^{-j}\delta \sqrt L) f \|_\infty \lesssim   \frac 1{\delta^{D/2}}\sum_{j\geq 0}  e^{- 2^{2j} \frac{ t}{8\delta^2}}
 2^{jD/2} 
 .$$

Put $A=  \frac{ t}{8\delta^2}$;
as:
$$
 \int_{2^j}^{2^{j+1}} x^{D/2} e^{-\frac A4x^2} \frac{Dx}x \geq 2^{\frac D2 j}  e^{-A2^{2j}} \log 2 $$
$$\sum_{j=0}^\infty 2^{\frac D2 j}  e^{-A2^{2j}} \leq \frac 1{\log 2} \int_{1}^{\infty} x^{D/2} e^{-\frac A4x^2} \frac{Dx}x 
=  \frac 1 {\log 2} \frac 12 (\frac 4A)^{D/4} \int_{A/4}^{\infty} u^{D/4} e^{-u} \frac{Du}u $$
as
$$ \hbox{ for all } a \in \R, X>0, \quad \int_X^\infty t^{a-1} e^{-t} dt \leq   2e^{-X} X^{a-1}, \quad \hbox{if} \;  X\geq 2 (a-1)$$
$$ \sum_{j=0}^\infty 2^{\frac d2 j}  e^{-A2^{2j}}  \leq \frac 4{\log 2} e^{-\frac A4} A^{-1}, \quad \hbox{if} \; A\geq 8( D-2)$$

So
$$\sum_{j\geq 0} \|\Psi( 2^{-j}\delta \sqrt L) f \|_\infty \leq   C  t^{-d/4}( \frac A4)^{d/4} e^{-\frac A4}( \frac A4)^{-1}. \qed$$ 

\paragraph{First step}: Fix $\delta$ such that $\| f-  \Phi( \delta \sqrt L)f \|_\infty <  \frac \epsilon2$ 
\\ 
Using the previous lemma, we need to choose $\delta $ so that  
\begin{equation}\label{aq1}
 \frac \epsilon 2 >   C  t^{-D/4}( \frac t{ 32\delta^{2}})^{D/4} e^{-\frac A4}( \frac A4)^{-1}, \quad \frac A4=  \frac{ t}{32 \delta^2}.
 \end{equation}
 Let us take :
 $$ \frac A4=\frac t{ 32\delta^2} = \alpha \log \frac 1\epsilon $$
then, as $ \epsilon^\mu \leq at.$
$$  C  t^{-D/4}( \frac A4)^{D/4} e^{-\frac A4}( \frac A4)^{-1}=   C  t^{-D/4} ( \alpha \log \frac 1\epsilon)^{D/4-1 }   \epsilon^{\alpha} 
\leq   C  (a\epsilon^\mu)^{-D/4} ( \alpha \log \frac 1\epsilon)^{D/4-1 }   \epsilon^{\alpha} \leq \frac \epsilon 2$$
if  $\alpha$ is suitably chosen. So for $\frac 1\delta \sim  \sqrt{\frac 1t \log \frac 1\epsilon}$,
$$\| f-  \Phi( \delta \sqrt L)f \|_\infty <  \frac \epsilon2$$

\paragraph{Second step }: $\epsilon-$ covering of $\bH_t^1$. 
\\
Now if $f \in \bH_t^1,$ using lemma \ref{lemm},  $  \|\Phi(\delta \sqrt L)f\|_\infty  \lesssim \frac{1}{t^{D/4}} $.
 Moreover  $\Phi(\delta \sqrt L)f \in \Sigma_{1/\delta},$
so, using (\ref{wav})
$$\Phi(\delta \sqrt L)f(x)=\sum_{\xi \in  \Lambda_{\gamma \delta}} \Phi(\delta \sqrt L)f(\xi) |B(\xi, \delta)| D^\delta_\xi(x) .$$
Let us consider the  following family  :
$$ f_{(k.)}=  C\sum_{\xi \in  \Lambda_{\gamma \delta}} k_\xi  \epsilon |B(\xi, \delta)| D^\delta_\xi(x)), \quad k_\xi \in \N, |k_\xi |\leq K \in \N, \quad  K C \epsilon \le \frac{1}{t^{D/4}}  $$
Certainly for all $ f \in \bH^1_t, $ there exists $( k_\xi)$  in the previous family such that 
\begin{align*} \|  \Phi(\delta \sqrt L)f -\sum_{\xi \in  \Lambda_{\gamma \delta}} k_\xi  C\frac \epsilon2 |B(\xi, \delta)| D^\delta_\xi(x)) \|_\infty
&=\| \sum_{\xi \in  \Lambda_{\gamma \delta}} (\Phi(\delta \sqrt L)f(\xi)- Ck_\xi \epsilon) |B(\xi, \delta)| D^\delta_\xi(x)\|_\infty
\\
&\lesssim  \sup_{\xi \in  \Lambda_{\gamma \delta}} |\Phi(\delta \sqrt L)f(\xi)- Ck_\xi \epsilon|   < \frac \epsilon2
\end{align*}
As $ \|\Phi(\delta \sqrt L)f  -f  \|_\infty \leq \frac \epsilon2$ , one can cover $\bH^1_t$ by balls centered in
 the $ f_{(k.)}$ of radius $\epsilon.$

The cardinality of this family of balls is : $(2K+1)^{card(\Lambda_{\gamma \delta})}$.
As $\gamma$ is a structural constant, $\epsilon^\mu \leq at$ and $\delta  \sim  \delta(t,\epsilon)$, clearly

$$ H(\epsilon, \bH^1_t, \L^\infty)  \lesssim     \cN(\delta(t,\epsilon), \M)  . \log \frac 1\epsilon $$

\subsection{ Bounds for $ \EE(\| W^t \|_\mb^2)$ }
In this section, we prove the following proposition with respect to the sup-norm. Similar  bounds in the $\L^2$ norm are obtained along the way (even slightly more precise).
\begin{proposition}
There exist universal constants $C_1$ and $C_2$ such that
\begin{equation}
C_1N(\sqrt{t},\M)\le \EE\|W^t \|_\infty^2\le C_2N(\sqrt{t},\M)\sup_{x\in \M}\frac 1{|B(x,\sqrt{t})|}
\label{w-infty}.
\end{equation}
\end{proposition}

We recall that $W^t$ writes
$$ W^t(x) = \sum_k \sum_{1\leq l \leq  \dim \cH_k} e^{-\lambda_kt/2}  X_k^l e^l_k(x)$$
where $X_k^l$ is a family of independent $N(0, 1)$ Gaussian variables. Clearly since $\M$ is supposed to have measure 1,
$$ \EE(\| W^t \|_2^2) \leq  \EE(\| W^t \|_\infty^2).$$
\\
As
$\| W^t \|_2^2 =  \sum_k  e^{-\lambda_kt}\sum_{1\leq l \leq  \dim \cH_k}   (X_k^l)^2 $, we get
$$\EE (\| W^t \|_2^2) =  \sum_k  e^{-\lambda_kt}  \dim \cH_k = Trace (e^{-tL}) = \int_{\M} P_t(u,u) d\mu(u).$$
Hence using Proposition \ref{prop5}
\begin{equation} \label{w-l2}
C'_1 2^{-2d}    N(\sqrt t, \M)  \leq  \EE (\| W^t \|_2^2)= \int_M P_t(u,u) d\mu(u)  \leq C_2' 2^{4d}  N(\sqrt t, \M).
\end{equation}
Now, let us first observe, using again  Proposition \ref{prop5}, that
\begin{align*}
\EE(\| W^t \|_\infty^2) &= \EE(\sup_{x\in \M} | W^t (x) |^2) 
\\
&\geq  \sup_{x\in \M} \EE( | W^t (x) |^2)\\
& =   \sup_{x\in \M} \EE(| \sum_k \sum_{1\leq l \leq  \dim \cH_k} e^{-\lambda_kt/2}  X_k^l e^l_k(x)|^2 
\\
&=\sup_{x\in \M} \sum_k \sum_{1\leq l \leq  \dim \cH_k} e^{-\lambda_kt}  ( e^l_k(x))^2
\\
&=
\sup_{x\in \M}  \sum_k  e^{-\lambda_kt} P_k(x,x)
\\
&=\sup_{x\in \M}P_t(x,x) \sim  \sup_{x\in \M} \frac 1{|B(x, \sqrt t)|}
\end{align*}
On the other side, using Cauchy-Schwarz inequality,
\begin{align*}
| W^t(x) |^2&= |\sum_k \sum_{1\leq l \leq  \dim \cH_k} e^{-\lambda_kt/2}  X_k^l e^l_k(x)|^2 
\\
&\leq \{ \sum_k \sum_{1\leq l \leq  \dim \cH_k} e^{-\lambda_k  t/2}  (X_k^l )^2\}
\{\sum_k \sum_{1\leq l \leq  \dim \cH_k} e^{-\lambda_k  t/2} (e^l_k(x))^2\}\\
&=\{ \sum_k \sum_{1\leq l \leq  \dim \cH_k} e^{-\lambda_k  t/2}  (X_k^l )^2\}
\{\sum_k  e^{-\lambda_k   t/2} P_k(x,x)\}
 \\
 &=\{ \sum_k \sum_{1\leq l \leq  \dim \cH_k} e^{-\lambda_k  t/2}  (X_k^l )^2\}
P_{  t/2} (x,x).
\end{align*}
So
\begin{align*}
\EE(\| W^t \|_\infty^2 ) &\leq  \EE  \{ \sum_k \sum_{1\leq l \leq  \dim \cH_k} e^{-\lambda_k  t/2}  (X_k^l )^2\} . \sup_{x \in \M} P_{  t/2} (x,x) 
\\
&= Trace (e^{-t/2 L}) \sup_{x \in \M} P_{  t/2} (x,x)\\
& = ( \int_{\M} P_{t/2}(u,u) d\mu(u) ) (\sup_{x \in \M} P_{  t/2} (x,x) ) .
\end{align*}
Hence, we get
 \begin{eqnarray*}
  \sup_{x \in \M} \frac 1{|B(x, \sqrt t)|} \sim \sup_{x\in \M}P_t(x,x)  \leq &  \EE(\| W^t \|_\infty^2) \leq  & ( \int_{\M} P_{t/2}(u,u) d\mu(u) ) (\sup_{x \in \M} P_{  t/2} (x,x) ) \\
  & & \sim    \cN(\sqrt{t}, \M) \sup_{x \in \M} \frac 1{|B(x, \sqrt t)|}. 
\end{eqnarray*}  
And we have in addition: $ N(\sqrt{t}, \M)  \sim \int_{\M}  \frac 1{|B(x, \sqrt t)|} d\mu(x)  \ll   \sup_{x \in \M} \frac 1{|B(x, \sqrt t)|}.$

\subsection{Lower bound for $A^t_f(\eps)$}

 \begin{theorem}\label{LowerboundA}
For $s >0$ fixed, there  exists $f \in B^s_{2,\infty}(\M), $ (the unit ball of the Besov space ) with $\|f \|_2^2 =1$ and constants $c>0,  C>0$ such that :
 $$\hbox{ for all } 1\geq  t>0, \; \hbox{ for all }  1>  \epsilon >0, \quad \inf_{\| f- h \|_2 \leq \epsilon} \| h \|_{\bH_t}^2 \geq C \epsilon^2 e^{ct\epsilon^{-2/s}}$$
 \end{theorem}

Let us take   $f$ such that $$  \|f \|_2 = 1  >  \epsilon >0.$$
%
We are interested in :
$$ \inf_{ \| f- P_{t/2} g\|_2= \epsilon } \|  g \|_2^2.$$
Let us put
\begin{equation}\label{1}
 \Phi(g) =\| f- P_{t/2} g\|^2_2=\| f \|_2^2 - 2\langle P_{t/2} f, g \rangle + \langle P_{t} g, g \rangle = \epsilon^2,\;  \Psi(g)= \| g \|^2_2.
 \end{equation}
 We have, 
$$D \Phi(g) = - 2P_{t/2} f + 2 P_t (g),\quad D\Psi(g) = 2g$$
%
$$\hbox{So, } \inf_{ \Phi(g) = \epsilon^2 }  \Psi(g)=  \Psi(g_0) \Longrightarrow  D\Psi(g_0)= -\mu D\Phi(g_0)$$
$$\hbox{ with } g_0 =  - \mu P_t(g_0) +\mu P_{t/2} f. $$
Necessarily $\mu \neq 0,$ otherwise $g_0=0$ and $\Phi(g_0)= \| f\|^2_2 \gg \epsilon^2.$
Let us put $\lambda = \frac 1\mu.$ We necessarily have
$ \lambda g_0 =   P_{t/2} f -  P_t(g_0)  $, hence
$( \lambda + P_{t}) (g_0) =    P_{t/2} f$, so
$$ g_0 = ( \lambda + P_{t})^{-1}  P_{t/2} f.$$
Let us now write the constraint :
$$\epsilon^2= \| f- P_{t/2} g\|^2_2= \| f- P_{t/2}  ( \lambda + P_{t})^{-1}  P_{t/2} f\|^2_2 = \| f-   ( \lambda + P_{t})^{-1}  P_t f\|^2_2
=\| \lambda  ( \lambda + P_{t})^{-1}   f\|^2_2  .$$
Clearly :
$$ \lambda \mapsto \| \lambda  ( \lambda + P_{t})^{-1}   f\|^2_2$$ is increasing from $0$ to $ \| f \|^2_2.$ As well,
$$ \lambda \mapsto 
   \|   ( \lambda + P_{t})^{-1}  P_{t/2} f\|^2_2$$
is decreasing .
On the other way : if $L = \int xdE_x$,  and 
$$\| \lambda  ( \lambda + P_{t})^{-1}   f\|^2_2 = \int_0^\infty ( \frac{\lambda}{\lambda + e^{-tx}})^2 d \langle E_x f ,f\rangle \geq \epsilon^2$$
and
$$  \| g_0\|_2^2 = \|   ( \lambda + P_{t})^{-1}  P_{t/2} f\|^2_2 =  \int_0^\infty ( \frac{1}{\lambda + e^{-tx}})^2 e^{-tx} d \langle E_x f ,f\rangle.$$

Let us recall the following result from \cite{CKP}, Lemma 3.19.
\begin{theorem}
There exists $ b >1, \; C"_1 >0 , \; C"_2 >0 , $ such that $ \hbox{ for all } \; \lambda \geq 1,
 \quad \delta = \frac 1\lambda, $  then
 $$( dim (\Sigma_{b\lam}) -dim (\Sigma_{\lam}) = dim (\Sigma_{b\lam} \ominus\Sigma_{\lam}) 
   =   \int_M  
P_{\Sigma_{b\lambda}}(x,x)  d\mu(x) -  \int_{\M}  
P_{\Sigma_\lambda}(x,x)  d\mu(x)  \neq 0$$
 and more precisely:
\begin{equation}\label{dimfinie}
  C"_1  \int_{\M}  \frac 1{|B(x, \delta)|} d\mu(x) \leq  dim (\Sigma_{b\lam} \ominus\Sigma_{\lam}) 
 \leq C"_2  \int_{\M}  \frac 1{|B(x, \delta)|} d\mu(x).
\end{equation}
\end{theorem}

As $ P_{\Sigma_{\sqrt a}}= E_a$, one can built a fonction $f \in \L^2$ such that :
$$  \| f - P_{\Sigma_{\sqrt a}}f\|_2^2=  \int_a^\infty \langle E_x f ,f\rangle = \| f \|_2^2 - \|E_a f  \|_2^2 = \| f - E_a f \|_2^2 = a^{-s}$$
for $ a= b^{2j},$  and $ j \in \N$.
It is enough to have :
$$  \| P_{\Sigma_{b^{j+1}} \ominus\Sigma_{b^j}} (f)\|^2_2 = b^{-2js}- b^{-2(j+1)s}$$ and this could be done by the previous theorem.

Let us choose for $\epsilon >0, \; b^{-2js} \geq  4\epsilon^2\geq  b^{-2(j+1)s}$.
So
$$ \int_{b^{2j}}^\infty \langle E_x f ,f\rangle =  b^{-2js} \geq 4 \epsilon^2\geq  b^{-2(j+1)s}= \int_{b^{2(j+1)}}^\infty \langle E_x f ,f\rangle$$
so, if $\lambda = e^{-ta}, \; a= b^{2j}$,
\begin{align*} \int_0^\infty ( \frac{\lambda}{\lambda + e^{-tx}})^2 d \langle E_x f ,f\rangle  &\geq 
 \int_a^\infty ( \frac{\lambda}{\lambda + e^{-tx}})^2 d \langle E_x f ,f\rangle \\
 &\ge 
 \int_a^\infty ( \frac{e^{-ta}}{e^{-ta} + e^{-tx}})^2 d \langle E_x f ,f\rangle 
 =
 \frac 14 \int_a^\infty  d \langle E_x f ,f\rangle \geq \epsilon^2.
 \end{align*}
But
\begin{align*}
 \| g_0 \|_2^2 &\geq \|   ( \lambda + P_{t})^{-1}  P_{t/2} f\|^2_2 =  \int_0^\infty ( \frac{1}{\lambda + e^{-tx}})^2 e^{-tx} d \langle E_x f ,f\rangle\\
 & =  e^{ta}  \int_0^\infty ( \frac{ 1}{e^{-ta} + e^{-tx}})^2 e^{-ta}e^{-tx} d \langle E_x f ,f\rangle =
 e^{ta}  \int_0^\infty ( \frac{ e^{-t/2x}  e^{-t/2a}}{e^{-ta} + e^{-tx}})^2  d \langle E_x f ,f\rangle \\
 &=  e^{ta}  \int_0^\infty ( \frac{1}{e^{-t/2(a-x)} + e^{-t/2(x-a)}})^2  d \langle E_x f ,f\rangle   
 \\
 &\geq e^{ta}  \frac 14 \int_0^\infty e^{-t |a-x|}  d \langle E_x f ,f\rangle  \geq e^{ta}  \frac 14 \int_{\frac a{b^2}}^a e^{-t (a-x)}  d \langle E_x f ,f\rangle  \geq  e^{t \frac a{b^2}} \frac 14 \int_{\frac a{b^2}}^a   d \langle E_x f ,f\rangle \\
 &=  e^{t \frac a{b^2}} \frac 14 \int_{b^{2j-2}}^{b^{2j}}   d \langle E_x f ,f\rangle = e^{t \frac a{b^2}} \frac 14 (b^{-(2j-2)s}-  b^{-2js})\\
 &=  e^{t \frac a{b^2}} \frac 14 b^{-2js} (b^{2s}-1) \geq \epsilon^2  (b^{2s}-1)  e^{t b^{2j-2}}
 \geq \epsilon^2  (b^{2s}-1) e^{t  c\epsilon^{-2/s}}; \quad c= 4^{-1/s} b^{-4}.
\end{align*} 
%
%
%


\section{ Appendix C: Compact Riemannian manifold.} \label{appendix-manifold}
 Let $\M$ be a compact Riemannian manifold without boundary. Let $\D(\M)$ be the algebra of infinitly differentiable functions. Associated to the Riemannian metric $\rho$, one defines a measure $dx$,
 a gradient operator $\nabla$ on $\D(\M)$ and the Laplace operator $\Delta$. It holds
 $$ \hbox{ for all } f,g \in \D(\M), \quad  \int_{\M} \Delta f(x) g(x) dx=- \int_{\M} |\nabla(f (x))|^2 dx.$$
 Thus $-\Delta$ is a positive symmetric operator. So actually
 $$ \L^2 = \oplus_{\lambda_k} \cH_{\lambda_k}, \quad \lambda_0=0 <\lambda_1< \lambda_2< ...$$
 $$  dim(\cH_{\lambda_k}) <\infty, \quad \cH_{\lambda_k} \subset \D(\M), \quad f\in  \cH_{\lambda_k}\Longleftrightarrow \Delta f=-\lambda_k f$$
 One can prove (see \cite{Grigoryan}) that the semi-group $e^{t\Delta}$ is a   kernel operator:
 $$ e^{t\Delta}(x,y)= \sum e^{-t\lambda_k } P_{\cH_{\lambda_k}}(x,y).$$
 Moreover this  kernel verifies:  there exist positive 
 $ C_1,C_2, C,c, \hbox{such that} $
$$ \hbox{ for all } u,v \in M, \;  \frac{C_2}{\sqrt{|B(u,\sqrt t)||B(v,\sqrt t)|}} e^{- C \frac{\rho^2(u,v)}t} \leq
 e^{tL} (u,v) \leq \frac{C_1}{\sqrt{|B(u,\sqrt t)||B(v,\sqrt t)|}} e^{- c \frac{\rho^2(u,v)}t}$$
 and moreover the property of doubling measure is verified. In fact we have a better result:
 \begin{proposition}
Let $\M$ be a compact Riemannian manifold of dimension $n.$ Then there exist $0<c\leq C <\infty$ such that :
$$\hbox{ for all } x\in \M, \hbox{ for all } 0<r<Diam(\M) , \quad cr^n \leq |B(x,r)| \leq Cr^n.$$
\end{proposition}
\noindent{\bf Proof :}

Let $\mu$ and $\rho$ be the (non normalized) Riemannian measure and metric on $ \M.$ The proposition is a consequence of the Bishop-Gromov comparison
 Theorem, see  \cite{Gromov} and \cite{Chavel}.
 
 As $\M$ is compact, clearly 
 $$\exists  \kappa \in \R, \quad \hbox{such that}: \hbox{ for all } x \in \M, \; 
 Ricc_x \geq (n-1) \kappa g_x$$
 where $Ricc$ is the Ricci tensor and $g$ is the metric tensor.
 Let $V_\kappa(r)$ be the volume of the (any) ball of radius $r$ in the model space of dimension $n$ and constant sectional curvature $\kappa.$
Let $V_n$ be  the volume of the unit ball of $\R^n.$
 \begin{enumerate}
\renewcommand\refname{Bibliography}
 \item For $\kappa >0, $ the model space is the sphere $\frac 1{\sqrt \kappa} \bS_n$ of $\R^{n+1}$ of radius $\frac 1{\sqrt \kappa}$ and
 $$ V_\kappa(r)= nV_n \int_0^r ( \frac{ \sin \sqrt{\kappa} t}{ \sqrt{\kappa} })^{n-1} dt; \;\hbox{so}\quad (\frac 2\pi)^{n-1} V_nr^n  \leq  V_\kappa(r) \leq V_nr^n$$
  \item For $\kappa =0, $ the model space is $\R^n$ and
$$ V_\kappa(r)= V_n r^n$$
 \item For $\kappa <0 $ the model space is the hyperbolic space of  constant sectional curvature $\kappa.$
$$ V_\kappa(r)= nV_n \int_0^r ( \frac{ \sinh \sqrt{|\kappa|} t}{ \sqrt{|\kappa|} })^{n-1} dt; \;\hbox{so}\quad  V_nr^n  \leq  V_\kappa(r) \leq V_nr^n e^{(n-1) \sqrt{|\kappa|} r}$$
 as $ s \leq \sinh(s) \leq se^s.$
\end{enumerate}
Moreover by the Bishop-Gromov comparaison comparison Theorem:
$ r \mapsto \frac{ |B(x,r)|}{V_\kappa(r)}$ is non increasing. So if $0<\epsilon <r  < s \leq R= diam(M) :$
$$\frac{ \mu(\M)}{V_\kappa(R)}=\frac{ |B(x,R)|}{V_\kappa(R)}   \leq \frac{ |B(x,s)|}{V_\kappa(s)} \leq \frac{ |B(x,r)|}{V_\kappa(r)} \leq \frac{ |B(x,\epsilon)|}{V_\kappa(\epsilon)}  \mapsto 1, \; \hbox{when} \; \epsilon \mapsto 0.$$
So
$$ \frac{V_\kappa(s)}{V_\kappa(r)}   \leq \frac{ |B(x,r)|}{ |B(x,s)|}; \quad \mu(M)\frac{  V_\kappa(r) }{V_\kappa(R)} \leq  |B(x,r)|\leq V_\kappa(r)$$
So
$$ A( \frac rs)^n \leq \frac{ |B(x,r)|}{ |B(x,s)|} \;\hbox{(doubling)} ;\quad  c r^n \leq   |B(x,r)| \leq  CV_n r^n, ; \; \hbox{(homogenity)}$$
$$  \hbox{for}\; \kappa > 0, C=1, \; c= (\frac 2\pi)^{n-1} \frac{\mu(\M)}{R^n}.   \quad A=(\frac 2\pi)^{n-1}.$$
$$ \hbox{for}\; \kappa = 0, C=1, \; c= \frac{\mu(\M)}{R^n}. \quad A=1.$$
$$\hbox{for}\; \kappa < 0, \; C= e^{(n-1) \sqrt{|\kappa|} R}; \; c= \frac{\mu(\M)}{R^n e^{(n-1) \sqrt{|\kappa|} R}} . \quad A= \frac{1}{ e^{(n-1) \sqrt{|\kappa|} R}}.$$

\begin{rem}
If $(\M, \mu, \rho)$ is a  compact metric space with a Borel measure $\mu$, then  if we have the doubling condition :
$$ 0<r<s \Longrightarrow  |B(x,s)| \leq  \frac 1A (\frac sr)^m |B(x,r)|$$ then 
$$ \hbox{ for all } r \leq R = diam(M), \;  Cr^m \leq |B(x,r)|, \quad \hbox{with} \quad C = \frac{A |\M|}{R^m}.$$
\end{rem}

\vspace{.5cm}

\noindent {\bf Acknowledgements.} The authors would like to thank Richard Nickl, Aad van der Vaart and Harry van Zanten for insightful comments on this work.

\bibliographystyle{abbrv}
\bibliography{bibsphere}

\end{document}